\documentclass[final,leqno]{siamltex}
\usepackage{amsmath,amssymb}
\usepackage{graphicx}
\newcommand{\bil}[1]{{#1}}

\newtheorem{remark_intro}{Remark}[section]
\newtheorem{remark_FD}{Remark}[section]
\newtheorem{remark}[remark_FD]{Remark}
\newtheorem{solver}{Algorithm}[section]

\newcommand{\hF}{\widehat{F}}
\newcommand{\cP}{\mathcal{P}}
\newcommand{\bP}{\mathbb{P}}
\newcommand{\hbP}{\mathbb{P}}
\newcommand{\cQ}{\mathcal{Q}}

\newcommand{\cS}{\mathcal{S}}
\newcommand{\cT}{\mathcal{T}}

\newcommand{\bfR}{\mathbf{R}}
\newcommand{\Ome}{\Omega}
\newcommand{\oOme}{\overline{\Omega}}
\newcommand{\ou}{\overline{u}}

\begin{document}

\title{Local discontinuous Galerkin methods for one-dimensional second order
fully nonlinear elliptic and parabolic equations}

\author{Xiaobing Feng\thanks{Department of Mathematics, The University
of Tennessee, Knoxville, TN 37996, U.S.A. (xfeng@math.utk.edu). The work of 
this author was partially supported by the NSF grant DMS-0710831.}
\and
Thomas Lewis\thanks{Department of Mathematics, The University
of Tennessee, Knoxville, TN 37996, U.S.A. (tlewis@math.utk.edu). 
The work of this author was partially supported by the NSF grant DMS-0710831.}
}

\maketitle
\date{\today}
\begin{abstract}
This paper is concerned with developing accurate and efficient 
discontinuous Galerkin methods for fully nonlinear second order  
elliptic and parabolic partial differential equations (PDEs) in the case
of one spatial dimension. The primary goal of the paper to develop
a general framework for constructing high order local discontinuous Galerkin (LDG) 
methods for approximating viscosity solutions of these fully nonlinear PDEs which 
are merely continuous functions by definition.
In order to capture discontinuities of the first order derivative $u_x$ 
of the solution $u$, two independent functions $q_1$ and  $q_2$ are introduced to 
approximate one-sided derivatives of $u$. Similarly, to capture the discontinuities 
of the second order derivative $u_{xx}$,  four independent functions $p_{1}$, $p_{2}$, $p_{3}$, 
and $p_{4}$ are used to approximate one-sided derivatives of $q_1$ and $q_2$.
The proposed LDG framework, which is based on a nonstandard mixed formulation of
the underlying PDE, embeds a given fully nonlinear problem into a mostly linear system 
of equations where the given nonlinear differential operator must be replaced by 
a numerical operator which allows multiple value inputs of the first and second order 
derivatives $u_x$ and $u_{xx}$. An easy to verify criterion for constructing 
``good" numerical operators is also proposed. It consists of a consistency and
a generalized monotonicity. To ensure such a generalized monotonicity, the crux of the 
construction is to introduce the numerical moment in the numerical operator, which plays a 
critical role in the proposed LDG framework. The generalized monotonicity gives the LDG methods
the ability to select the viscosity solution among all possible solutions.
The proposed framework extends a companion finite difference framework
developed by the authors in \cite{FKL11} and allows for the approximation of 
fully nonlinear PDEs using high order polynomials and non-uniform meshes. 
Numerical experiment results are also presented to demonstrate 
the accuracy, efficiency and utility of the proposed LDG methods. 
\end{abstract}

\begin{keywords}
Fully nonlinear PDEs, viscosity solutions, local discontinuous Galerkin methods,
\end{keywords}

\begin{AMS}
65N30, 
65M60, 
35J60, 
35K55, 
\end{AMS}

\pagestyle{myheadings}
\thispagestyle{plain}
\markboth{XIAOBING FENG AND THOMAS LEWIS}{LDG METHODS FOR
SECOND ORDER FULLY NONLINEAR PDEs}

\section{Introduction}\label{sec-1}

This is the third paper in a series \cite{FKL11,Feng_Lewis12b} which is 
devoted to developing finite difference (FD) and discontinuous 
Galerkin (DG) methods for approximating {\em viscosity solutions}
of the following general one-dimensional fully nonlinear second order elliptic 
and parabolic equations:
\begin{equation}\label{pde_ell}
F \left( u_{x x}, u_x, u, x \right) = 0 , \qquad x \in \Omega := (a,b) \subset \bfR,
\end{equation}
and
\begin{equation} \label{pde}
u_t + F \left( u_{x x}, u_x, u, x,t \right) = 0 , 
\qquad (x,t) \in \Omega_T:=\Omega\times (0,T),
\end{equation}
which are complemented by appropriate boundary and initial conditions. 

Fully nonlinear PDEs, which are nonlinear in the highest order derivatives of 
the solution functions in the equations,  arise in many applications such as 
antenna design, astrophysics, differential geometry, fluid mechanics, image processing, 
meteorology, mesh generation, optimal control, optimal mass transport, etc 
(cf.  \cite[section 5]{FGN12}), and,
as a result, the solution of each of these application problems critically depends
on the solution of their underlying fully nonlinear PDEs. In particular, it calls for 
efficient and reliable numerical methods for computing the viscosity solutions of these fully 
nonlinear PDEs. Currently, the availability of such numerical methods are very limited
(cf. \cite{FGN12}).

The goal of this paper is to design and implement a class of local discontinuous 
Galerkin (LDG) methods for the fully nonlinear equations \eqref{pde_ell} and 
\eqref{pde}.  The more involved high dimensional generalizations of the LDG methods 
of this paper will be reported in \cite{Feng_Lewis12c}.

Because of the full nonlinearity, integration by parts, which is the necessary tool for
constructing any DG method, cannot be performed on equation \eqref{pde_ell}.
{\em The first key idea} of this paper is to introduce the auxiliary variables $p:=u_{x x}$ and $q:= u_x$ and 
rewrite the original fully nonlinear PDE in the following nonstandard mixed form:  
\begin{eqnarray} \label{mixed_1}
F(p,q,u,x) &=0, \\
q- u_x &=0, \label{mixed_2} \\
p-q_x &=0. \label{mixed_3}
\end{eqnarray}

Unfortunately, since $u_x$ and $u_{x x}$ may not exist for a viscosity 
solution $u\in C^0(\Omega)$, the the above mixed form may 
not make sense. To overcome this difficulty, {\em our second key idea} 
is to replace $q=u_x$ by two possible values of 
$u_x$, namely, its left and right limits, and $p= q_x$ by two possible
values for each possible $q$.  Thus, we have the auxiliary variables
$q_1, q_2 : \Omega \to \bfR$ and $p_1, p_2, p_3, p_4 : \Omega \to \bfR$
such that
\begin{eqnarray} 
q_1(x)-u_x(x^-) &=0, \label{qdn} \\
q_2(x)-u_x(x^+) &=0, \label{qup} \\
p_1(x)-q_{1x}(x^-) &=0, \label{pdndn} \\
p_2(x)-q_{1x}(x^+) &=0, \label{pdnup} \\
p_3(x)-q_{2x}(x^-) &=0, \label{pupdn} \\
p_4(x)-q_{2x}(x^+) &=0. \label{pupup} 
\end{eqnarray}
We remark that \eqref{qdn} paired with the equation \eqref{pdndn} or \eqref{pdnup}, 
and \eqref{qup} paired with equation \eqref{pupdn} or \eqref{pupup}, can each be regarded 
as a ``one-sided" Poisson problem in $u$ (in a mixed form) with source terms
$p_1$, $p_2$, $p_3$, $p_4$, respectively.

To incorporate the multiple values of $u_x$ and $u_{x x}$, 
equation \eqref{mixed_1} must be modified because $F$ is only defined for 
single value functions $p$ and $q$. To this end, we need 
{\em the third key idea} of this paper, that is, to replace \eqref{mixed_1} by 
\begin{eqnarray} \label{F_hat}
\hF(p_1, p_2, p_3, p_4, q_1, q_2, u, x) = 0,
\end{eqnarray} 
where $\hF$, which is called a {\em numerical operator}, should be some
well-chosen approximation to $F$. 

Natural questions now arise regarding to the choice of $\hF$. 
For example, what are criterions for $\hF$ and
how to construct such a numerical operator? 
These are two immediate questions which must be 
addressed. To do so, we need {\em the fourth key idea} of this paper, 
which is to borrow and adapt the notion of the numerical 
operators from our previous work \cite{FKL11} where a general
finite difference framework has been developed for fully nonlinear 
second order PDEs. In summary, the criterions for $\hF$ include 
{\em consistency} and {\em g-monotonicity} (generalized monotonicity), 
for which precise definitions can be found in section \ref{sec-2}.  
It should be pointed out that in order to construct the desired numerical operator
$\hF$, a fundamental idea used in \cite{FKL11} is to introduce the 
concept of {\em the numerical moment}, which can be regarded as 
a direct numerical realization for the moment term in {\em the vanishing
moment methodology} introduced in \cite{Feng_Neilan08} (also see
\cite[section 4]{FGN12},\cite{Feng_Neilan11}). 

Finally, we need to design 
a DG discretization for the mixed system \eqref{qdn}--\eqref{F_hat} 
to complete the construction of our LDG method.  This then calls for 
{\em the fifth key idea} of this paper, which is to use different {\em numerical fluxes} 
in the formulations of LDG methods for the four ``one-sided" Poisson 
problems in their mixed forms described by \eqref{qdn} -- \eqref{pupup}. 
We remark that, to the best of our knowledge, this is one of a few scenarios 
in numerical PDEs where the flexibility and superiority (over other numerical 
methodologies) of the DG methodology makes a vital difference.
 
This paper consists of four additional sections. In section \ref{sec-2}
we collect some preliminaries including the definition of viscosity 
solutions, the definitions of the consistency and g-monotonicity of 
numerical operators, and the concept of the numerical moment.  
In section \ref{sec-3} we give the detailed formulation of LDG methods 
for fully nonlinear elliptic equation \eqref{pde_ell} following the outline 
described above. In section \ref{sec-4} we consider both explicit and 
implicit in time fully discrete LDG methods for fully nonlinear 
parabolic equation \eqref{pde}. The explicit four stage classical Ronge-Kutta method
and the implicit trapezoidal method combined with the spatial LDG methods will be 
explicitly constructed. In section \ref{sec-5} we present a number of numerical 
experiments for the proposed LDG methods for the fully nonlinear 
elliptic equation \eqref{pde_ell} and their fully discrete 
counterparts for the parabolic equation \eqref{pde}. 
These numerical experiments not only verify the accuracy of the proposed 
LDG methods but also demonstrate the efficiency and utility of these methods.

\section{Preliminaries} \label{sec-2}
Standard space notations are adopted in this paper. For example,  $B(\Ome)$, 
$USC(\Ome)$ and $LSC(\Ome)$ denote, respectively, the spaces of bounded,
upper semi-continuous, and lower semicontinuous functions on $\Ome$.
For any $v\in B(\Ome)$, we define
\[
v^*(x):=\limsup_{y\to x} v(y) \qquad\mbox{and}\qquad
v_*(x):=\liminf_{y\to x} v(y). 
\]
Then, $v^*\in USC(\Ome)$ and $v_*\in LSC(\Ome)$, and they are called
{\em the upper and lower semicontinuous envelopes} of $v$, respectively.
If $v$ is continuous,  there obviously holds $v=v^*=v_*$. 

Let  $F: \cS^{d\times d}\times\mathbf{R}^d\times \mathbf{R}\times \oOme \to \mathbf{R}$ be a 
bounded function, where $\cS^{d\times d}$ denotes the set of $d\times d$ symmetric real matrices. 
Then, the general second order fully nonlinear PDE takes the form
\begin{align}\label{e2.1}
F(D^2u,\nabla u, u, x) = 0 \qquad\mbox{in } \oOme.
\end{align}
Note that here we have used the convention of writing the boundary condition as a
discontinuity of the PDE (cf. \cite[p.274]{Barles_Souganidis91}).

The following two definitions can be found in \cite{Gilbarg_Trudinger01,
Caffarelli_Cabre95,Barles_Souganidis91}.

\begin{definition}\label{def2.1}
Equation \eqref{e2.1} is said to be elliptic if, for all
$(\mathbf{q},\lambda,x)\in \mathbf{R}^d\times \mathbf{R}\times \oOme$, there holds
\begin{align}\label{e2.2}
F(A, \mathbf{q}, \lambda, x) \leq F(B, \mathbf{q}, \lambda, x) \qquad\forall 
A,B\in \cS^{d\times d},\, A\geq B, 
\end{align}
where $A\geq B$ means that $A-B$ is a nonnegative definite matrix.
\end{definition}
We note that when $F$ is differentiable, the ellipticity
also can be defined by requiring that the matrix $\frac{\partial F}{\partial A}$
is negative semi-definite (cf. \cite[p. 441]{Gilbarg_Trudinger01}).

\begin{definition}\label{def2.2}
A function $u\in B(\Ome)$ is called a viscosity subsolution (resp.
supersolution) of \eqref{e2.1} if, for all $\varphi\in C^2(\oOme)$,
if $u^*-\varphi$ (resp. $u_*-\varphi$) has a local maximum
(resp. minimum) at $x_0\in \oOme$, then we have
\[
F_*(D^2\varphi(x_0),\nabla \varphi(x_0), u^*(x_0), x_0) \leq 0 
\]
(resp. $F^*(D^2\varphi(x_0),\nabla \varphi(x_0), u_*(x_0), x_0) \geq 0$).
The function $u$ is said to be a viscosity solution of \eqref{e2.1}
if it is simultaneously a viscosity subsolution and a viscosity
supersolution of \eqref{e2.1}.
\end{definition}

We note that if $F$ and $u$ are continuous, then the upper and lower $*$
indices can be removed in Definition \ref{def2.2}. The definition
of ellipticity implies that the differential operator $F$
must be non-increasing in its first argument in order to be
elliptic. It turns out that ellipticity provides a sufficient 
condition for equation \eqref{e2.1} to fulfill a maximum principle
(cf. \cite{Gilbarg_Trudinger01, Caffarelli_Cabre95}).
From the above definition it is clear that viscosity solutions
in general do not satisfy the underlying PDEs in a tangible sense, and
the concept of viscosity solutions is {\em nonvariational}. Such
a solution is not defined through integration by parts against arbitrary test
functions; hence, it does not satisfy an integral identity. This nonvariational 
nature of viscosity solutions is the main obstacle that prevents direct construction
of Galerkin-type methods, which require variational formulations to start. 

\smallskip
The following definitions are adapted from \cite{FKL11} in the case $d=1$.

\smallskip
\begin{definition}\label{def2.3}
\begin{itemize}
\item[{\rm (i)}] A function $\hF: \bfR^8\to \bfR$ is called  
a {\em numerical operator}. 
\item[{\rm (ii)}] A numerical operator $\hF$ is said to be consistent (with 
the differential operator $F$) if $\hF$ satisfies
\begin{align}\label{A1a}
\liminf_{p_k\to p, k=1,2,3,4\atop q_1,q_2 \to q, \lambda_1\to \lambda,\xi_1\to \xi} 
\hF(p_1,p_2,p_3,p_4,q_1,q_2,\lambda_1, \xi_1) \geq F_*(p,q,\lambda,\xi),\\
\limsup_{p_k\to p, k=1,2,3,4 \atop q_1,q_2 \to q, \lambda_1\to \lambda,\xi_1\to \xi} 
\hF(p_1,p_2,p_3,p_4, q_1,q_2, \lambda_1,\xi_1) 
\leq F^*(p,q,\lambda,\xi), \label{A1b}
\end{align}
where $F_*$ and $F^*$ denote respectively the lower and the upper
semi-continuous envelopes of $F$.

\item[{\rm (iii)}] A numerical operator $\hF$ is said to be
{\em g-monotone} if $\hF(p_1,p_2,p_3,p_4,q_1,q_2,\lambda,\xi)$ is monotone increasing 
in $p_1$ and $p_4$ and monotone decreasing in $p_2$ and $p_3$, that is, 
$\hF(\uparrow,\downarrow,\downarrow,\uparrow,q_1,q_2,\lambda,\xi)$.
\end{itemize}

\end{definition}

\medskip
We remark that the above consistency and g-monotonicity play a critical role in 
the finite difference framework established in \cite{FKL11}. They also play
an equally critical role in the LDG framework of this paper. 
In practice, the consistency is easy to fulfill and 
to verify, but the g-monotonicity is not. In order to ensure the 
g-monotonicity, one key idea of \cite{FKL11} is to introduce 
 {\em the numerical moment} to help.  The following is an example of 
a so-called Lax-Friedrichs-like numerical operator adapted from \cite{FKL11}:
\begin{align}\label{LF1}
\hF(p_1, p_2, p_3, p_4 ,q_1, q_2,\lambda, \xi)
&:= F(\frac{p_2 + p_3}{2},\frac{q_1 + q_2}{2},\lambda,\xi) \\
\nonumber & \qquad\quad
+ \alpha \bigl(p_1 - p_2 - p_3 + p_4 \bigr),
\end{align}
where $\alpha \in \bfR$ is an undetermined positive constant and 
the last term in \eqref{LF1} is called {\em the numerical moment}. 
It is trivial to verify that $\hF$ is consistent with $F$.
By choosing $\alpha$ correctly, we can also ensure g-monotonicity. For example, 
suppose $F$ is differentiable and there exists a positive constant $M$ such that
\begin{equation} \label{alpha_bound} 
M > \left| \frac{\partial F}{\partial u_{x x}} \right|. 
\end{equation}
Then, it is trivial to check that the Lax-Friedrichs-like numerical operator 
is g-monotone provided that  $\alpha\geq M$. 

We conclude this section with a few remarks about the above definitions.

\smallskip
\begin{remark_intro}
(a) By the definition of the ellipticity for $F$, the monotonicity constraints on $\hF$ with respect
to $p_2$ and $p_3$ in the definition of g-monotonicity are natural.

(b) By choosing the numerical moment correctly, the numerical operator $\hF$ then behaves 
like a uniformly elliptic operator, even if the PDE operator $F$ is a degenerate elliptic operator.
The consistency assumption then guarantees that the numerical operator is still a reasonable
approximation for the PDE operator.

(c) Sometimes it may not be feasible to globally bound $\frac{\partial F}{\partial u_{x x}}$; 
however, it is sufficient to chose a value for $\alpha$ such that the g-monotonicity property is 
preserved locally over each iteration of the nonlinear solver for a given initial guess.

(d) The role of the numerical moment as well as the interpretation of the numerical moment
will be further discussed in Section~\ref{alpha_tests}.
\end{remark_intro}

\section{Formulation of LDG methods for elliptic problems}\label{sec-3}
 
We first consider the elliptic problem $(\ref{pde_ell})$ with boundary conditions
\begin{equation}
u(a) = u_a \quad\mbox{and}\quad  u(b) = u_b \label{bc_ell} 
\end{equation}
for two given constants $u_a$ and $u_b$. 

Let $\left\{ x_j \right\}_{j=0}^J \subset \overline{\Omega}$ be a mesh 
for $\overline{\Omega}$ such that $x_0 = a$ and $x_J = b$.  
Define $I_{j} = \left( x_{j-1} , x_{j} \right)$ and $h_j = x_{j} - x_{j-1}$
for all $j = 1, 2, \ldots, J$, $h_0 = h_{J+1} = 0$, and 
$h = \max_{1 \leq j \leq J} h_j$.  Let $\mathcal{T}_h$ denote
the collection of the intervals $\{I_j\}_{j=1}^J$ which form a partition 
of the domain $\overline{\Omega}$. We also introduce the broken $H^1$-space 
and broken $C^0$-space
\[
H^1(\cT_h):= \prod_{I\in \cT_h} H^1(I), \qquad C^0(\cT_h):=\prod_{I\in \cT_h} C^0(\overline{I}),
\]
and the broken $L^2$-inner product
\[
(v ,w)_{\mathcal{T}_h}:= \sum_{j=1}^J \int_{I_j} v w\, dx \qquad 
\forall v,w\in L^2(\Ome).
\]
For a fixed integer $r \geq 0$, we define the standard DG finite element space
$V^h \subset H^1(\cT_h)\subset L^2 (\mathcal{T}_h)$ as 
\[
V^h := \prod_{I \in \mathcal{T}_h} \mathcal{P}_{r} (I),
\]
where $ \mathcal{P}_{r} (I)$ denotes the set of all polynomials on $I$ with
degree not exceeding $r$.  We also introduce the following
standard jump notation:
\begin{align*}
[v_h(x_j)] &:= v_h(x_j^-)-v_h(x_j^+) \qquad\mbox{for } j=1,2,\cdots, J-1.
\end{align*}

We now are ready to formulate our LDG discretizations for equations 
\eqref{qdn}--\eqref{F_hat}. First, for (fully) nonlinear equation \eqref{F_hat}
we simply approximate it by its broken $L^2$-projection into $V^h$, namely,
\begin{equation}\label{pde_ell_weak}
\bil{a}_0 \bigl(u_h , q_{1h}, q_{2h}, p_{1h}, p_{2h}, p_{3h}, p_{4h}; \phi_{0h} \bigr)  = 0 
\qquad \forall \phi_{0h} \in V^h, 
\end{equation}
where
\[
\bil{a}_0 (u_h, q_{1h}, q_{2h}, p_{1h}, p_{2h}, p_{3h}, p_{4h}; \phi_{0h}) 
=\Bigl(\hF(p_{1h}, p_{2h}, p_{3h}, p_{4h}, q_{1h}, q_{2h}, u_h,\cdot),\phi_{0h} \Bigr)_{\mathcal{T}_h}. 
\]

Next, we discretize the four {\em linear} equations \eqref{pdndn}--\eqref{pupup}.
Notice that for given ``sources" $\{p_i\}_{i=1}^4$, \eqref{qdn} and \eqref{pdndn}, 
\eqref{qdn} and \eqref{pdnup}, \eqref{qup} and \eqref{pupdn}, and \eqref{qup} and \eqref{pupup} 
are four (different) Poisson equations for $u$. Thus, we can use the mixed upwinding LDG formulation 
for the Laplacian operator to discretize these equations. The only difference in the four equations will be 
how we choose our upwinding numerical fluxes for $u_h$, $q_{1h}$ and $q_{2h}$.
To realize the above strategy, we first define the element-wise LDG formulation, and we then
define the whole domain LDG formulation afterward. 

\subsection{Element-wise LDG formulation} \label{local_DG}

Suppose that values for $u_h(a^-)$, $u_h(b^+)$, $q_{ih}(a^-)$, and $q_{ih}(b^+)$ for $i=1,2$
are given.  
We postpone explaining how these values are chosen until section~\ref{boundary_fluxes}.  
Our LDG discretization of equations \eqref{qdn}--\eqref{pupup} is defined as follows:  for all $\phi_{ih}\in V^h$, 
\begin{equation} \label{q_local}
\int_{I_j} q_{ih} \, \phi_{ih} \, dx + \int_{I_j} u_h \, (\phi_{ih})_x \, dx
= u_h(x_j^\sigma) \, \phi_{ih}(x_j^-) - u_h(x_{j-1}^\sigma) \, \phi_{ih} (x_{j-1}^+)
\end{equation}
for $i = 1,2$, $j=1,2,\cdots,J$, and 
\[
\sigma = \begin{cases}
	- & \text{if } i = 1, \\
	+ & \text{if } i = 2;
\end{cases}
\] 
and for all $\psi_{kh} \in V^h$, 
\begin{equation} \label{p_local}
\int_{I_j} p_{kh} \, \psi_{kh} \, dx + \int_{I_j} q_{\hat{k}h} \, (\psi_{kh})_x\, dx
= q_{\hat{k}h}(x_j^\beta) \, \psi_{kh}(x_j^-) - q_{\hat{k}h}(x_{j-1}^\beta) \, \psi_{kh} (x_{j-1}^+)
\end{equation}
for $k = 1,2,3,4$, and 
\[
\hat{k} = \begin{cases}
	1 & \text{if } k \text{ is odd}, \\
	2 & \text{if } k \text{ is even},
\end{cases} 
\qquad 
\beta = \begin{cases}
	- & \text{if } \hat{k} = 1, \\
	+ & \text{if } \hat{k} = 2. 
\end{cases}
\]
Notice that the nodal values determined by $\sigma$ and $\beta$ follow directly from equations $(\ref{qdn}) - (\ref{pupup})$.

\subsection{Boundary numerical fluxes} \label{boundary_fluxes}

To complete the construction, we must specify how the boundary numerical 
flux values for $u_h$, $q_{1h}$, and $q_{2h}$ are determined in the above formulation.  
Due to the inherent jumps of piecewise constant functions, which corresponds to 
the case $r = 0$,  we shall consider the two cases  $r \geq 1$ and $r = 0$ separately. 

When $r \geq 1$, we have freedom to control how the functions $u_h$, 
$q_{1h}$, and $q_{2h}$ approach the boundary.  
Thus, we can assume continuity across the boundary for $u_h$, $q_{1h}$, and $q_{2h}$.
Considering the boundary conditions given by \eqref{bc_ell}, the continuity requirement 
naturally leads to 
\begin{equation} \label{ubc_r1}
	u_h(a^\pm) = u_a, \qquad u_h(b^\pm) = u_b.
\end{equation}

On the other hand, since no boundary data for $q_{1h}$ or $q_{2h}$ is given,  any choice of 
the boundary numerical fluxes for them is a guess (unless one already knows 
the exact solution $u$).  Here we choose
\begin{equation} \label{qbc_r1}
	q_{ih}(a^-) = q_{ih}(a^+), \qquad q_{ih}(b^+) = q_{ih}(b^-),  \qquad i=1,2. 
\end{equation}
It is important to note that  both $q_{ih}(a^+)$ and $q_{ih}(b^-)$ ($i=1,2$) are 
treated as unknowns in the above LDG formulation. The choice \eqref{qbc_r1} is equivalent to
requiring that $q_{ih}$ is continuous at the boundary nodes $x=a$ and $x=b$.

We now consider the case $r=0$. To define the boundary numerical fluxes, we first examine the 
consequences of the interior flux choices represented by the above LDG formulation.
Suppose $\mathcal{T}_h$ is a uniform mesh and denote the midpoint of $I_j$ by $\hat{x}_j$ 
for all $I_j \in \mathcal{T}_h$. Define $U_j := u^h \left( \hat{x}_{j} \right)$.
Then,  it follows from \eqref{q_local} and \eqref{p_local} that 
\begin{align} 
	q_{1h} \left( \hat{x}_j \right) & = \frac{ U_j - U_{j-1} }{h} := \delta_x^- U_j, \label{q1h_r0} \\
	q_{2h} \left( \hat{x}_j \right) & = \frac{ U_{j+1} - U_j }{h} := \delta_x^+ U_j,   \label{q2h_r0} \\
	p_{1h} \left( \hat{x}_j \right) & = \frac{q_{1h}(\hat{x}_{j}) - q_{1h}(\hat{x}_{j-1})}{h}
		= \frac{ U_{j-2} - 2 U_{j-1} + U_j }{h^2} 
		:= \delta_x^2 U_{j-1} , \label{p1h_r0} \\
	p_{2h} \left( \hat{x}_j \right) & = \frac{q_{1h}(\hat{x}_{j+1}) - q_{1h}(\hat{x}_{j})}{h}
		= \frac{ U_{j-1} - 2 U_j + U_{j+1} }{h^2}
		:= \delta_x^2 U_{j} , \\
	p_{3h} \left( \hat{x}_j \right) & = \frac{q_{2h}(\hat{x}_{j}) - q_{2h}(\hat{x}_{j-1})}{h}
		= \frac{ U_{j-1} - 2 U_j + U_{j+1} }{h^2}
		:= \delta_x^2 U_{j} , \\
	p_{4h} \left( \hat{x}_j \right) & = \frac{q_{2h}(\hat{x}_{j+1}) - q_{2h}(\hat{x}_{j})}{h}
		= \frac{ U_{j} - 2 U_{j+1} + U_{j+2} }{h^2}
		:= \delta_x^2 U_{j+1} , \label{p4h_r0}
\end{align}
for $j = 3, 4, \ldots, J-2$. Thus, in order to define boundary values for $u_h$, $q_{1h}$, and $q_{2h}$, 
we need to define ghost values $U_{-1}$, $U_0$, $U_{J+1}$, and $U_{J+2}$ that are equivalent to 
extending the solution $u$ to the outside of the domain $\Omega$. Below we describe a natural way 
to do such an extension that is consistent with the interpretation of the auxiliary variables.

From the Dirichlet boundary data for $u$, a natural choice is that $U_{0} = u_a$ and $U_{J+1} = u_b$.
This is equivalent to assuming 
\begin{equation} \label{ubc_r0}
u^h(a^-) = u_a , \qquad u^h(b^+) = u_b.
\end{equation}
In other words, we extend the boundary data for $u$ away from the boundary over an interval 
of length $h$. Due to the inherent discontinuities of the piecewise constant functions,  
$u^h(a^+)$ and $u^h(b^-)$ are treated as unknowns. Otherwise, the boundary data would 
be extended into the interior of the domain over an interval of length $h$.

From \eqref{q1h_r0} and \eqref{q2h_r0} we see that $q_{2h} \left( \hat{x}_j \right) = q_{1h} \left( \hat{x}_{j+1} \right)$
and $q_{1h} \left( \hat{x}_j \right) = q_{2h} \left( \hat{x}_{j-1} \right)$ in the interior of the domain.   
Extending this relationship to the boundary yields
\begin{align}  \label{qbc_r0a}
	q_{2h} (a^-) = q_{1h} (a^+) , \qquad
	q_{1h} ( b^+) = q_{2h} (b^-) ,
\end{align}
where both $q_{1h}(a^+)$ and $q_{2h}(b^-)$ are treated as unknowns in the above LDG formulation. 

Finally, we need to define values for $q_{1h}(a^-)$ and $q_{2h}(b^+)$.  
Using \cite{Feng_Neilan08} as a guide, we are led to choosing
\begin{equation} \label{qbc_r0c}
q_{1h}(a^-) = q_{1h}(a^+) , \qquad q_{2h}(b^+) = q_{2h}(b^-).
\end{equation}
We note that this is consistent with discretizing the auxiliary boundary conditions 
\begin{equation}\label{qbc_r0d}
(q_{1h})_x(a) = (q_{2h})_x(b) = 0.
\end{equation} 
In order words,  we require that $q_{1h}$ and $q_{2h}$ are constant across the boundary.
Using ghost values, the above requirements are equivalent to imposing the constraints
\begin{equation*}
U_{-1} = 2 u_a - U_1, \qquad U_{J+2} = 2 u_b - U_J.
\end{equation*}

\begin{remark_FD}
From the imposed boundary conditions we can see that the relationship
$p_{2h} = p_{3h}$ has been extended to the boundaries.
Thus, using the ghost values defined above and substituting the equations 
\eqref{q1h_r0}--\eqref{p4h_r0} into \eqref{pde_ell_weak}, 
we successfully recover the convergent finite difference method defined in \cite{FKL11} for
the grid function $U$. Thus, for $r=0$, the convergence of the proposed LDG method 
is obtained.  Heuristically, using  higher order elements should increase the rate and/or accuracy 
of convergence.
\end{remark_FD}

\subsection{Whole domain LDG formulation} \label{global_DG}

Using the above element-wise LDG formulation  \eqref{q_local} and \eqref{p_local}, 
and substituting the boundary numerical flux values from section~\ref{boundary_fluxes}, 
we get the following whole domain LDG discretization of \eqref{qdn}--\eqref{pupup}:
\begin{align}
\left( q_{ih} , \phi_{ih} \right)_{\cT_h} + \bil{a}_i (u_h, \phi_{ih} ) 
 &= \bil{f}_i (\phi_{ih}), 
&& \forall \phi_{ih} \in V^h, \,\, i=1,2, \label{q_bil} \\
\left( p_{jh} , \psi_{jh} \right)_{\cT_h} + \bil{b}_j (q_{1h}, q_{2h};\psi_{jh} ) 
& = 0, 
&& \forall \psi_{jh} \in V^h, \,\, j=1,2,3,4, \label{p_bil}
\end{align}
where 
\begin{align*}
&\bil{a}_1 (v, \varphi)
=( v , \varphi_x)_{\cT_h} - (1-\kappa_r) \, v(b^-) \, \varphi(b^-)
-\sum_{j=1}^{J-1} v (x_{j}^-) \bigl[ \varphi (x_{j}) \bigr] , \\
&\bil{a}_2 (v, \varphi)
=( v , \varphi_x)_{\cT_h} + (1-\kappa_r) \, v(a^+) \, \varphi(a^+)
-\sum_{j=1}^{J-1} v (x_{j}^+) \bigl[ \varphi (x_{j}) \bigr] , \\
&\bil{b}_1 (v_1, v_2; \varphi)
=( v_1 , \varphi_x)_{\cT_h} + v_1(a^+) \, \varphi(a^+) - v_1(b^-) \, \varphi(b^-)
-\sum_{j=1}^{J-1} v_1 (x_{j}^-) \bigl[ \varphi (x_{j}) \bigr] , \\
&\bil{b}_4 (v_1, v_2; \varphi)
=( v_2 , \varphi_x)_{\cT_h} + v_2(a^+) \, \varphi(a^+) - v_2(b^-) \, \varphi(b^-)
-\sum_{j=1}^{J-1} v_2 (x_{j}^+) \bigl[ \varphi (x_{j}) \bigr] , \\
&\bil{b}_2 (v_1, v_2; \varphi)
=( v_1 , \varphi_x)_{\cT_h} + v_1(a^+) \, \varphi(a^+) 
- (1 - \kappa_r) \, v_2(b^-) \, \varphi(b^-) \\
&\hskip 1in 
- \kappa_r \, v_1(b^-) \, \varphi(b^-)  -\sum_{j=1}^{J-1} v_1 (x_{j}^+) \bigl[ \varphi (x_{j}) \bigr] , \\
&\bil{b}_3 (v_1, v_2; \varphi) =( v_2 , \varphi_x)_{\cT_h} 
+ (1 - \kappa_r) \, v_1(a^+) \, \varphi(a^+)  
+ \kappa_r \, v_2(a^+) \, \varphi(a^+)   \\
&\hskip 1in 
 - v_2(b^-) \, \varphi(b^-)  -\sum_{j=1}^{J-1} v_2 (x_{j}^-) \bigl[ \varphi (x_{j}) \bigr] , 
 \end{align*}
 and
 \begin{align*}
\bil{f}_1 (\phi) & = \kappa_r \, u_b \, \phi(b^-) - u_a \, \phi(a^+) , \\
\bil{f}_2 (\phi) & =  u_b \, \phi(b^-) - \kappa_r \, u_a \, \phi(a^+) , 
\end{align*}
for
\begin{equation}\label{kapr}
\kappa_r = \begin{cases}
	0 & \text{if } r=0, \\
	1 & \text{otherwise} .
\end{cases}
\end{equation}
 

In summary, our nonstandard LDG methods for the fully nonlinear Dirichlet problem 
\eqref{pde_ell} and \eqref{bc_ell} are defined as seeking
$\bigl(u_h , q_{1h}, q_{2h}, p_{1h}, p_{2h}, p_{3h}, p_{4h} \bigr) 
\in \left( V^h \right)^7$
such that \eqref{pde_ell_weak}, \eqref{q_bil}, and \eqref{p_bil} hold.  
 
We conclude the section with a few remarks.

\smallskip
\begin{remark}
(a) Looking backwards, \eqref{q_bil} and \eqref{p_bil} provide a proper interpretation for 
each of $q_{ih}$ and $p_{jh}$ for $i=1,2$ and $j = 1,2,3,4$, for a given function $u_h$. Each 
$q_{ih}$ defines a discrete derivative for $u_h$ and each $p_{jh}$ defines a 
discrete second order derivative for $u_h$.
The functions $q_{1h}$ and $q_{2h}$ should be very close to each other if $u_x$ exists.
Similarly, the functions $p_{1h}$, $p_{2h}$, $p_{3h}$, and $p_{4h}$ 
should be very close to each other if $u_{x x}$ exists.
However, their discrepancies are expected to be large 
if $u_x$ or $u_{x x}$, respectively, do not exist. 
The auxiliary functions $q_{ih}$ defined by \eqref{q_bil} and the auxiliary functions
$p_{jh}$ defined by \eqref{p_bil}
can be regarded as high order extensions of their lower order finite difference counterparts 
defined in \cite{FKL11}.

(b) It is easy to check that the linear equations defined by \eqref{q_bil}--\eqref{p_bil}
are linearly independent.

(c) Notice that \eqref{pde_ell_weak}, \eqref{q_bil}, and \eqref{p_bil} form a nonlinear 
system of equations, with the nonlinearity only appearing in $\bil{a}_0$. 
Thus, a nonlinear solver is necessary in implementing the above scheme. 
In section \ref{sec-5}, an iterative method is used with initial guess given by the 
linear interpolant of the boundary data.  
Since a good initial guess is essential for most nonlinear solvers 
to converge, another possibility is to first linearize the nonlinear 
operator and solve the resulting linear system first.  
However, we show in our numerical tests that the simple initial 
guess works well in many cases. We suspect that the g-monotonicity 
of $\hF$ enlarges the domain of ``good" initial guesses over which 
the iterative method converges.
\end{remark}

\section{Formulation of fully discrete LDG methods for parabolic problems}
\label{sec-4}

The goal of this section is to extend the LDG methods for elliptic problems to solving the
initial-boundary value problem \eqref{pde} using the method of lines.  Let the initial condition be given by
\begin{equation}\label{ic}
u(0,x) = u_0(x), \qquad \forall x \in \Omega,
\end{equation}
and the boundary conditions be given by
\begin{equation} \label{bc_time}
u(t,a) = u_a(t), \qquad u(t,b) = u_b(t), \qquad \forall t \in (0,T].
\end{equation}
We shall consider both the implicit trapezoidal rule and the (explicit) fourth order classical Runge-Kutta method (i.e., RK4) 
for the time-discretization. In practice, the time-discretization scheme should be chosen to match the order 
of the spacial discretization.  Thus, when using a piecewise constant element, a sufficient choice for the time
discretization would be forward or backward Euler.  However, when using  higher order elements, 
a higher order scheme such as RK4 should be chosen.

We first formulate the semi-discrete in space discretization for equation \eqref{pde},
which is a straightforward adaptation of the one described in section \ref{sec-3}.  Let $\phi \in V^h$.  
Replacing the PDE operator with a numerical operator in \eqref{pde}, 
and using the LDG framework of section \ref{sec-3}, we obtain the semi-discrete equation
\begin{equation}\label{semi-discrete}
\bigl( u_{ht}, \phi_h \bigr)_{\cT_h} = 
 - \bigl( \hF \left( p_{1h}, p_{2h}, p_{3h}, p_{4h}, q_{1h}, q_{2h}, u_h, t, \cdot \right) , 
		\phi_h \bigr)_{\cT_h},  \quad \forall \phi_h \in V^h , 
\end{equation}
where, given $u_h$ at time $t$, corresponding values for $q_{ih}$ and $p_{jh}$ can be found using 
the methodology below.

We now describe a full discretization procedure for \eqref{pde} by applying an ODE solver
to the semi-discrete equations given in \eqref{semi-discrete}.
For a fixed integer $M>0$,  let $\Delta t = \frac{T}{M}$ and  $t_k := k \, \Delta t$ for $k \in (0,M]$.
Notationally, $u_h^n(x) \in V^h$ will be an approximation for $u(t_n, x)$ for $n = 0, 1, \ldots, M$.
We define the initial value $u_h^0$ to be the $L^2$-projection of $u_0$,  namely,  
\begin{equation} \label{proj_ic}
\bigl( u_h^0 , \phi_h \bigr)_{\cT_h} = \bigl( u_0, \phi_h \bigr)_{\cT_h}, 
	\qquad \forall \phi_h \in V^h.
\end{equation}

Next, we introduce several ``one-sided" discrete differential operators, 
which will be used to define explicit time-stepping schemes and to define 
the auxiliary variables at time $t=0$ in implicit time-stepping schemes. 
We first define two ``one-sided" first order discrete derivatives
$\cQ_1^k v, \cQ_2^k v \in V^h$ for a given function $v(\cdot, t_k)\in H^1(\cT_h)$ by
\begin{align} 
\bigl(\cQ_1^k v , \phi_{1} \bigr)_{\cT_h} 
	& = \kappa_r \, u_b(t_k) \, \phi_{1}(b^-) - u_a(t_k) \, \phi_{1}(a^+) 
		+ (1-\kappa_r) \, v(b^-) \, \phi_{1}(b^-) \label{q1_time} \\
	\nonumber & \qquad - \bigl( v, \phi_{1x} \bigr)_{\cT_h}
		+ \sum_{j=1}^{J-1} v(x_j^-) \, \bigl[ \phi_{1}(x_j) \bigr] , 
		\qquad \forall \phi_{1} \in V^h , \\
\bigl(\cQ_2^k v , \phi_{2} \bigr)_{\cT_h} 
	& = u_b(t_k) \, \phi_{2}(b^-) - \kappa_r \, u_a(t_k) \, \phi_{2}(a^+) 
		- (1-\kappa_r) \, v(a^+) \, \phi_{2}(a^+) \label{q2_time} \\
	\nonumber & \qquad - \bigl( v, \phi_{2x} \bigr)_{\cT_h}
		+ \sum_{j=1}^{J-1} v(x_j^+) \, \bigl[ \phi_{2}(x_j) \bigr] , 
		\qquad \forall \phi_{2} \in V^h.
\end{align}
The above definitions are inspired  by \eqref{q_bil}. The super-index $k$ on $Q_i^k$  indicates that
the definitions are $t$-dependent because of the boundary terms. 

 We also define four discrete ``one-sided" second order discrete derivatives 
 $\cP_1^k v, \cP_2^k v,$  $\cP_3^k v , \cP_4^k v \in V^h$ at time $t_k$ by
\begin{align}
\bigl(\cP_1^k v , \psi_{1} \bigr)_{\cT_h} 
	& = \cQ_1^k v(b^-) \, \psi_{1}(b^-) - \cQ_1^k v(a^+) \, \psi_{1}(a^+) \label{p1_time} \\
	\nonumber & \qquad - \bigl( \cQ_1^k v, \psi_{1x} \bigr)_{\cT_h}
		+ \sum_{j=1}^{J-1} \cQ_1^k v(x_j^-) \, \bigl[ \psi_{1}(x_j) \bigr] , 
		\qquad \forall \psi_{1} \in V^h, \\
\bigl(\cP_4^k v , \psi_{4} \bigr)_{\cT_h} 
	& = \cQ_2^k v(b^-) \, \psi_{4}(b^-) - \cQ_2^k v(a^+) \, \psi_{4}(a^+) \label{p4_time} \\
	\nonumber & \qquad - \bigl( \cQ_2^k v, \psi_{4x} \bigr)_{\cT_h}
		+ \sum_{j=1}^{J-1} \cQ_2^k v(x_j^+) \, \bigl[ \psi_{4}(x_j) \bigr] , 
		\qquad \forall \psi_{4} \in V^h,  \\
\bigl(\cP_2^k v , \psi_{2} \bigr)_{\cT_h} 
	& = (1-\kappa_r) \, \cQ_2^k v(b^-) \, \psi_{2}(b^-) + \kappa_r \, \cQ_1^k v(b^-) \, \psi_{2}(b^-)  \\
	\nonumber & \qquad
		- \cQ_1^k v(a^+) \, \psi_{2}(a^+) \label{p2_time} 
	 - \bigl( \cQ_1^k v, \psi_{2x} \bigr)_{\cT_h} \\
	 \nonumber & \qquad
		+ \sum_{j=1}^{J-1} \cQ_1^k v(x_j^+) \, \bigl[ \psi_{2}(x_j) \bigr] , 
		\qquad \forall \psi_{2} \in V^h,  \\ 
\bigl(\cP_3^k v , \psi_{3} \bigr)_{\cT_h} 	& = \cQ_2^k v(b^-) \, \psi_{3}(b^-) 
		- (1-\kappa_r) \, \cQ_1^k v(a^+) \, \psi_{3}(a^+)  \\
         \nonumber & \qquad
		- \kappa_r \, \cQ_2^k v(a^+) \, \psi_{3}(a^+)  \label{p3_time} 
	       - \bigl( \cQ_2^k v, \psi_{3x} \bigr)_{\cT_h} \\
         \nonumber & \qquad
 		+ \sum_{j=1}^{J-1} \cQ_2^k v(x_j^-) \, \bigl[ \psi_{3}(x_j) \bigr] , 
		\qquad \forall \psi_{3} \in V^h,
\end{align}
where $\kappa_r$ is defined by \eqref{kapr}.  The above four definitions are motivated by \eqref{p_bil}. 

Lastly, to simplify the presentation, we introduce the operator notation
\begin{equation}\label{fhat_short}
\hF_k [v]  
:= \hF \left( \cP_1^k v, \cP_2^k v, \cP_3^k v , \cP_4^k v, \cQ_1^k v , \cQ_2^k v , v , t_k , x \right) .
\end{equation}
Using the new notation, the semi-discrete equation can be rewritten compactly as 
\begin{equation}\label{semi}
\bigl( u_{ht}(t_k, \cdot), \phi_h \bigr)_{\cT_h} = 
 - \bigl( \hF_k \left[ u_h(t_k, \cdot) \right] , \phi_h\bigr)_{\cT_h} ,  \quad \forall \phi_h\in V^h, \,
\, k=1,2,\cdots, M.
\end{equation}

\subsection{The fourth order classical Runge-Kutta method} \label{RK4}

A straightforward application of the fourth order classical Runge-Kutta (RK4) method to \eqref{semi} yields 
\[
\bigl( u_h^n, \phi_h \bigr)_{\cT_h} = 
\bigl( u_h^{n-1}, \phi_h \bigr)_{\cT_h} 
	+ \frac{1}{6} \bigl( \xi_1 + 2 \xi_2 + 2 \xi_3 + \xi_4 , \phi_h \bigr)_{\cT_h} , 
\quad \forall \phi \in V^h, \, n=1,2,\cdots, M,
\]
where
\begin{align*}
\bigl( \xi_1 , \phi_h\bigr)_{\cT_h} & = - \Delta t \bigl( \hF_{n-1} [ u_h^{n-1} ] , \phi_h \bigr)_{\cT_h} , \\ 
\bigl( \xi_2 , \phi_h \bigr)_{\cT_h} & = 
	- \Delta t \bigl( \hF_{n-\frac12} [ u_h^{n-1} + \frac12 \xi_1] , \phi_h \bigr)_{\cT_h} , \\ 
\bigl( \xi_3 , \phi_h \bigr)_{\cT_h} & = 
	- \Delta t \bigl( \hF_{n-\frac12} [  u_h^{n-1} + \frac12 \xi_2 ] , \phi_h \bigr)_{\cT_h} , \\ 
\bigl( \xi_4 , \phi_h \bigr)_{\cT_h} & = 
	- \Delta t \bigl( \hF_{n} [  u_h^{n-1} + \xi_3 ] , \phi_h \bigr)_{\cT_h} .
\end{align*}

Notice that in the above explicit time-stepping scheme, the function $u_h^{n}$ 
is defined as an $L^2$-projection of the source data based on $u_h^{n-1}$. 
However, the boundary conditions are not enforced in the definition for $u_h^n$. 
To take care of the boundary conditions, we choose to enforce them  
weakly,  which requires the introduction of a modified $L^2$-projection. 
Specifically, for any $v\in L^2(\Ome)$, we recall that the standard 
$L^2$-projection $\bP_h v\in V^h$ of $v$ is defined by 
\begin{equation} \label{L2_proj}
\bigl( \bP_h v, \phi_h \bigr)_{\cT_h} = \bigl(v, \phi_h \bigr)_{\cT_h} , 
\qquad\forall \phi_h \in V^h.
\end{equation}
For any $v\in C^0(\cT_h)$, we introduce a modified $L^2$-projection 
$\hbP_h^k : L^2(\Ome)\cap C^0(\cT_h)\to V^h$ at time $t_k \in (0, T]$
by
\begin{align} \label{mL2_proj}
&\bigl( \hbP_h^k v, \phi_h \bigr)_{\cT_h} 
+ \frac{1}{\sqrt{h}} \Bigl(\hbP_h^k v(a) \phi_h(a) + \hbP_h^k v (b) \phi_h(b) \Bigr) \\
&\hskip 0.5in
= \bigl(v, \phi_h \bigr)_{\cT_h} + \frac{1}{\sqrt{h}} \Bigl(u_a(t_k) \, \phi_h(a^+) + u_b(t_k) \, \phi_h(b^-) \Bigr) , 
\qquad\forall \phi_h\in V^h. \nonumber
\end{align}
Clearly, the boundary conditions \eqref{bc_time} are weakly enforced in \eqref{mL2_proj} 
via a penalty technique, an idea which dates back to Nitsche \cite{Nitsche70}.  

Using the above discrete differential operators $\cQ_i^k$, $\cP_j^k$, $\bP_h^k$, and $\bP_h$
for $i=1,2$ and $j=1,2,3,4$, and using the notation given in \eqref{fhat_short},  
our fully-discrete RK4 method for the initial-boundary value problem \eqref{pde}, \eqref{bc_time}, 
and \eqref{ic} is defined as follows:
for $n = 1,2, \ldots, M$, 
\begin{align} \label{rk4}
u_h^{n} & = \hbP_h^{n} \Bigl( u_h^{n-1} + \frac16( \xi_1^{n-1} + 2 \xi_2^{n-1} + 
	2 \xi_3^{n-1} + \xi_4^{n-1} ) \Bigr) , \\
\xi_1^{n-1} & = - \Delta t \,  \bP_h  \hF_{n-1} [ u_h^{n-1}] , \\ 
\xi_2^{n-1} & = - \Delta t \, \bP_h  \hF_{n-\frac12} [ u_h^{n-1} + \frac12\xi_1^{n-1}  ] , \\ 
\xi_3^{n-1} & = - \Delta t \, \bP_h  \hF_{n-\frac12} [  u_h^{n-1} + \frac12 \xi_2^{n-1} ]  , \\ 
\xi_4^{n-1} & = - \Delta t \, \bP_h  \hF_{n} [  u_h^{n-1} + \xi_3^{n-1} ]  , \\
u_h^0 & = \bP_h u_0 .
\end{align}
We remark that it is easy to verify that the value $\xi_4^{n-1}$ actually already take into account the 
boundary conditions at time $t_n$ because the evaluation calls $Q_1^n$ and $Q_2^n$, which  
in turn have the boundary conditions built-in.  Thus, in the above formulation, the boundary condition 
enforcement can actually be successfully relaxed by replacing $\hbP_h^n$ with $\bP_h$ in \eqref{rk4}.  
However, for other explicit methods  a weak boundary condition enforcement
method such as the above modified $L^2$-projection is necessary, 
especially if the boundary conditions are not consistent with the initial condition. For example,
in the forward Euler method, defined by 
\begin{equation}\label{fEuluer}
u_h^n=\hbP_h^n \bigl( u_h^{n-1} -\Delta t \hF_{n-1}[u_h^{n-1}]\bigr), 
\end{equation}
we can see that the approximation at time $t_n$ relies upon the modified $L^2$-projection in order 
to see the Dirichlet boundary condition at the current time.

\subsection{The trapezoidal method}
\label{trapezoidal}

Applying the trapezoidal rule to \eqref{semi},  we obtain
\[
\Bigl( u_h^n + \frac{\Delta t}{2}  \hF_n [ u_h^n ], \phi \Bigr)_{\cT_h} = 
\Bigl( u_h^{n-1} -\frac{\Delta t}{2} \hF_{n-1} [ u_h^{n-1} ], \phi \Bigr)_{\cT_h} , 
\qquad \forall \phi \in V^h ,  
\]
and  $n = 1, 2, \ldots , M$. Thus, using the trapezoidal rule to discretize \eqref{semi}, 
and using the implicit equalities
\[
q_{1h}^n = \cQ_1^n q_{1h}^n , \qquad
q_{2h}^n = \cQ_2^n q_{2h}^n , 
\]
and
\[
p_{1h}^n = \cP_1^n p_{1h}^n , \qquad
p_{2h}^n = \cP_2^n p_{2h}^n , \qquad
p_{3h}^n = \cP_3^n p_{3h}^n , \qquad
p_{4h}^n = \cP_4^n p_{4h}^n ,
\]
for $n = 1, 2, \ldots, M$, 
the fully discrete trapezoidal LDG method for approximating solutions to 
\eqref{pde}, \eqref{ic}, and \eqref{bc_time} is defined by seeking
$(u_h^n, q_{1h}^n, q_{2h}^n, p_{1h}^n, p_{2h}^n, p_{3h}^n, p_{4h}^n) \in (V^h)^7$
such that 
\begin{align}
&\Bigl( u_h^n + \frac{\Delta t}{2}  \hF_n [ u_h^n ], \phi_{0h} \Bigr)_{\cT_h} = 
\Bigl( u_h^{n-1} -\frac{\Delta t}{2} \hF_{n-1} [ u_h^{n-1} ], \phi_{0h} \Bigr)_{\cT_h} , 
\quad \forall \phi \in V^h,\\
&\bigl( q_{ih}^n , \phi_{ih} \bigr)_{\cT_h} 
	+ \widehat{a}_i \left( u_h^n , \phi_{ih} \right) = g_i \left( t^n, \phi_{ih} \right) ,
	\quad \forall \phi_{ih} \in V^h, i = 1, 2, \label{ai_time} \\
&\bigl( p_{jh}^n , \psi_{jh} \bigr)_{\cT_h} 
	+ \widehat{b}_j \left( q_{1h}^n, q_{2h}^n ; \psi_{jh} \right) = 0 ,
	\quad \forall \psi_{jh} \in V^h, j = 1, 2, 3, 4 \label{bj_time} , 
\end{align}
where $u^0_h = \bP_h u_0$, $q_{ih}^0 = \cQ_i^0 u_h^0$ for $i=1,2$, 
$p_{jh}^0 = \cP_j^0 u_h^0$ for $j = 1,2,3,4$,  and 
\begin{align*}
&\widehat{a}_1 (v^n, \varphi)
=( v^n , \varphi_x)_{\cT_h} - (1-\kappa_r) \, v^n(b^-) \, \varphi(b^-)
-\sum_{j=1}^{J-1} v^n (x_{j}^-) \bigl[ \varphi (x_{j}) \bigr] , \\
&\widehat{a}_2 (v^n, \varphi)
=( v^n , \varphi_x)_{\cT_h} + (1-\kappa_r) \, v^n(a^+) \, \varphi(a^+)
-\sum_{j=1}^{J-1} v^n (x_{j}^+) \bigl[ \varphi (x_{j}) \bigr] , \\
&\widehat{b}_1 (v_1^n, v_2^n; \varphi)
=( v_1^n , \varphi_x)_{\cT_h} + v_1^n(a^+) \, \varphi(a^+) - v_1^n(b^-) \, \varphi(b^-)
-\sum_{j=1}^{J-1} v_1^n (x_{j}^-) \bigl[ \varphi (x_{j}) \bigr] , \\
&\widehat{b}_4 (v_1^n, v_2^n; \varphi)
=( v_2^n , \varphi_x)_{\cT_h} + v_2^n(a^+) \, \varphi(a^+) - v_2^n(b^-) \, \varphi(b^-)
-\sum_{j=1}^{J-1} v_2^n (x_{j}^+) \bigl[ \varphi (x_{j}) \bigr] , \\
&\widehat{b}_2 (v_1^n, v_2^n; \varphi)
=( v_1^n , \varphi_x)_{\cT_h} + v_1^n(a^+) \, \varphi(a^+) 
- (1 - \kappa_r) \, v_2^n(b^-) \, \varphi(b^-)  \\
&\hskip 1.5in 
- \kappa_r \, v_1^n(b^-) \, \varphi(b^-) -\sum_{j=1}^{J-1} v_1^n (x_{j}^+) \bigl[ \varphi (x_{j}) \bigr] , \\
&\widehat{b}_3 (v_1^n, v_2^n; \varphi) =( v_2^n , \varphi_x)_{\cT_h} 
+ (1 - \kappa_r) \, v_1^n(a^+) \, \varphi(a^+) + \kappa_r \, v_2^n(a^+) \, \varphi(a^+) \\
&\hskip 1.5in 
- v_2^n(b^-) \, \varphi(b^-) -\sum_{j=1}^{J-1} v_2^n (x_{j}^-) \bigl[ \varphi (x_{j}) \bigr] , \\
&g_1 (t^n, \phi)  = \kappa_r \, u_b(t^n) \, \phi(b^-) - u_a(t^n) \, \phi(a^+) , \\
&g_2 (t^n, \phi)  =  u_b(t^n) \, \phi(b^-) - \kappa_r \, u_a(t^n) \, \phi(a^+) .
\end{align*}
Again,  $\kappa_r$ is defined by \eqref{kapr}.  Notice that the above fully discrete formulation 
amounts to approximating a nonhomogeneous fully nonlinear elliptic equation at each time step 
using the LDG method defined in section~\ref{sec-3}.

\section{Numerical experiments}\label{sec-5}

In this section, we present a series of numerical tests to demonstrate 
the utility of the proposed LDG methods for  fully nonlinear PDEs of
the types \eqref{pde_ell} and $(\ref{pde})$. In all of our tests we shall use 
uniform spatial meshes as well as uniform temporal meshes for the time-dependent 
problems. To solve the resulting nonlinear algebraic systems, we use 
the Matlab built-in nonlinear solver {\em fsolve} for the job.
For the elliptic problems we choose the initial guess as the 
linear interpolant of the boundary data. 
For parabolic problems, we let $u^0_h = \bP_h u_0$,  and then define all auxiliary
variables by $q_{ih}^0 = \cQ_i^0 u_h^0$ for $i=1,2$ and 
$p_{jh}^0 = \cP_j^0 u_h^0$ for $j = 1,2,3,4$. 
 Also, the initial guess for {\em fsolve} at the $n$th time step will be chosen as 
 the computed solution at the previous time step when using implicit methods.
The role of the numerical moment will be further explored in section~\ref{alpha_tests}.

For our numerical tests, errors will be measured in the $L^\infty$ norm 
and the $L^2$ norm, where the errors are measured at the current time 
step for the time-dependent problems. For both elliptic and parabolic
test problems where the error is dominated by the spatial discretization errors,  
it appears that the spatial errors are of order $\mathcal{O} (h^{r+1})$.
Thus, the schemes appear to exhibit an optimal rate of convergence in both norms.

\subsection{Elliptic test problems}
\label{ell_tests_1d}
We first present the results for four test problems of type \eqref{pde_ell}.
Both Monge-Amp\`ere and Bellman types of equations will be tested. 

\smallskip
{\bf Test 1.} 
Consider the elliptic Monge-Amp\`{e}re problem
\begin{align*}
- u_{x x}^2 + 1 & = 0 , \qquad 0 < x < 1 , \\
u(0) = 0, \quad  u(1) & = \frac12 .
\end{align*}
It is easy to check that this problem has exactly two classical solutions: 
\[
u^+ (x) = \frac{1}{2} x^2 , \qquad u^- (x) = - \frac{1}{2} x^2 + x , 
\]
where $u^+$ is convex and $u^-$ is concave.  Note that $u^+$ is the unique 
viscosity solution which we want our numerical schemes to converge to.
In section~\ref{alpha_tests} we shall give some insights about 
the selectiveness of our schemes. 

\begin{figure}
\centerline{
\includegraphics[scale=0.35]{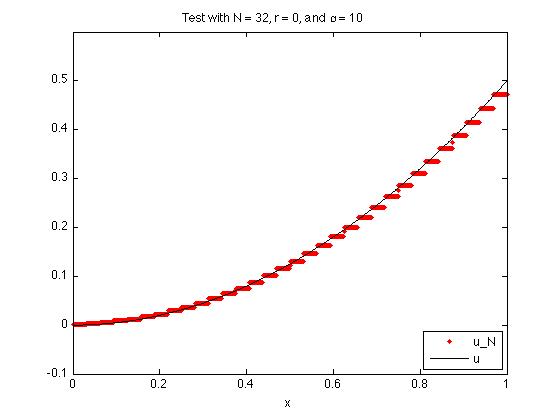}
\includegraphics[scale=0.35]{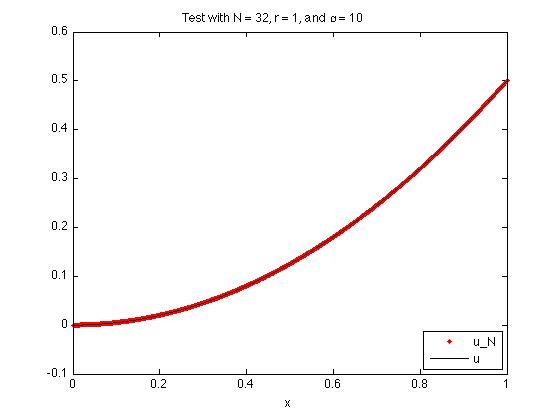}
}
\centerline{
{\small
\begin{tabular}{| c | c | c | c c | c c | c c |} 
		\hline
	$r$ & Norm & $h = 1/4$ 
		&\multicolumn{2}{| c |}{$h = 1/8$}
		&\multicolumn{2}{| c |}{$h = 1/16$}
		&\multicolumn{2}{| c |}{$h = 1/32$} \\ 
		\cline{3-9} 
	&   & Error & Error & Order & Error & Order & Error & Order  \\ 
		\hline \cline{1-9}
	0 & $L^{2}$ & 7.1e-02 & 3.5e-02 & 1.02 & 1.4e-02 
	  & 1.30 & 7.5e-03 & 0.92  \\ 
	  & $L^{\infty}$ & 1.3e-01 & 8.7e-02 & 0.57 & 5.3e-02 
	  & 0.73 & 2.9e-02 & 0.87  \\ 
		\hline
	1 & $L^{2}$ & 1.6e-02 & 5.0e-03 & 1.67 & 1.3e-03 
	  & 1.90 & 3.4e-04 & 1.95  \\ 
	  & $L^{\infty}$ & 2.2e-02 & 6.3e-03 & 1.84 & 1.6e-03 
	  & 2.00 & 3.9e-04 & 2.00  \\ 
		\hline
	2 & $L^{2}$ & 3.1e-13 & 3.0e-13 & 0.03 & 3.0e-13 
	  & -0.01 & 3.1e-13 & -0.01  \\ 
	  & $L^{\infty}$ & 7.4e-13 & 6.1e-13 & 0.28 & 6.7e-13 
	  & -0.14 & 7.1e-13 & -0.08  \\    
		\hline
\end{tabular}
}
}
\caption{Test 1 with $\alpha = 10$.}
\label{ma_1d_r}
\end{figure}

We approximate the given problem for various degree elements  ($r = 0,1,2$) 
to see how the approximation converges with respect to $h$.  Note, when $r=0,1$, the 
solution is not in the DG space $V^h$. The numerical results 
are shown in Figure~\ref{ma_1d_r}.  We observe that
the approximations using $r=2$ are almost exact for each mesh size.  
This is expected since $u^+ \in V^h$ when $r = 2$.

\begin{figure} 
\centerline{
\includegraphics[scale=0.35]{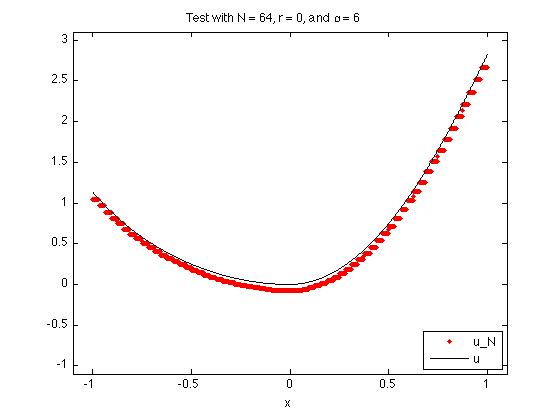}
\includegraphics[scale=0.35]{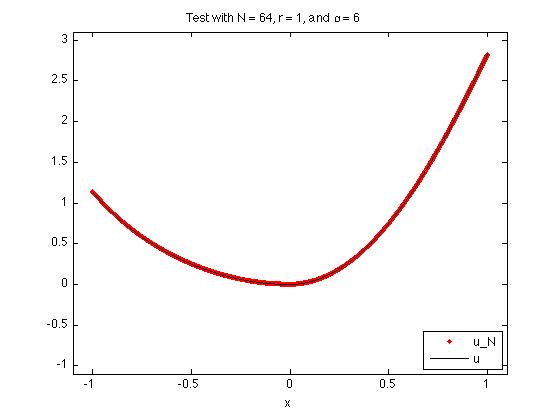}
}
\centerline{
{\tiny
\begin{tabular}{| c | c | c | c c | c c | c c | c c |} 
		\hline
	$r$ & Norm & $h = 1/2$ 
		&\multicolumn{2}{| c |}{$h = 1/4$}
		&\multicolumn{2}{| c |}{$h = 1/8$}
		&\multicolumn{2}{| c |}{$h = 1/16$}
		&\multicolumn{2}{| c |}{$h = 1/32$} \\ 
		\cline{3-11} 
	&   & Error & Error & Order & Error & Order & Error & Order & Error & Order  \\ 
		\hline \cline{1-11}
	0 & $L^{2}$ & 1.8 & 9.0e-01 & 1.04 & 4.3e-01 
	  & 1.05 & 2.1e-01 & 1.04 
	  & 1.0e-01 & 1.02  \\ 
	  & $L^{\infty}$ & 2.5 & 1.2 & 0.99 & 6.2e-01 
	  & 1.00 & 3.1e-01 & 1.00 
	  & 1.6e-01 & 1.00  \\ 
		\hline
	1 & $L^{2}$ & 2.9e-01 & 6.3e-02 & 2.20 & 1.9e-02 
	  & 1.73 & 7.0e-03 & 1.44 
	  & 2.8e-03 & 1.34  \\ 
	  & $L^{\infty}$ & 2.9e-01 & 6.4e-02 & 2.17 & 2.0e-02 
	  & 1.66 & 7.6e-03 & 1.42 
	  & 3.1e-03 & 1.30  \\ 
		\hline
	2 & $L^{2}$ & 5.7e-03 & 8.2e-04 & 2.80 & 1.3e-04 
	  & 2.66 & 3.2e-05 & 2.03 
	  & 9.1e-06 & 1.81  \\ 
	  & $L^{\infty}$ & 2.0e-02 & 3.1e-03 & 2.70 & 4.2e-04 
	  & 2.87 & 5.5e-05 & 2.94 
	  & 8.0e-06 & 2.77  \\ 
		\hline
	3 & $L^{2}$ & 8.8e-04 & 7.7e-05 & 3.51 & 3.0e-06 
	  & 4.68 & 1.4e-07 & 4.42 
	  & 1.0e-08 & 3.76  \\ 
	  & $L^{\infty}$ & 2.1e-03 & 1.4e-04 & 3.90 & 8.6e-06 
	  & 4.01 & 5.6e-07 & 3.94 
	  & 9.5e-08 & 2.57  \\ 
		\hline
\end{tabular}
}
}
\caption{Test 2 with $\alpha = 6$.}
\label{test2_1d}
\end{figure}

\medskip
{\bf Test 2.} 
Consider the problem
\begin{align*}
	- u_{ x x}^3 + u_{x x} + S(x)^3 - S(x) & = 0, \qquad -1 < x < 1 , \\
	u(-1) = - \sin (1) - 8 \cos(0.5) + 9 , & \quad
	u(1) = \sin (1) - 8 \cos(0.5) + 9 ,     
\end{align*}
where
\[
	S(x) = \begin{cases}
			\frac{2 x}{|x|} \cos(x^2) - 4 x^2 \sin (x |x|) + 2 \cos(\frac{x}2) + 2 , & x \neq 0 , \\
			- 4 x^2 \sin (x |x|) + 2 \cos(\frac{x}2) + 2 , & x = 0 .
		   \end{cases}
\]
This problem has the exact solution $u(x) = \sin \left( x | x | \right) - 8 \cos \left( \frac{x}2 \right) + x^2 + 8$.  
Note that this problem is not monotone decreasing in $u_{x x}$, and the exact solution is not
twice differentiable at $x = 0$.  
However, the derivative of $F$ with respect to $u_{x x}$ is uniformly bounded.
The numerical results are shown in Figure~\ref{test2_1d}.
As expected, we can see from the plot that the error appears largest 
around the point $x = 0$, and both the accuracy and order of 
convergence improve as the order of the element increases.  
For finer meshes, we see the rates of convergence begin to deteriorate.  
Theoretically, we expect less than optimal rates of convergence due to the low 
regularity of the solution.

\begin{figure} 
\centerline{
\includegraphics[scale=0.35]{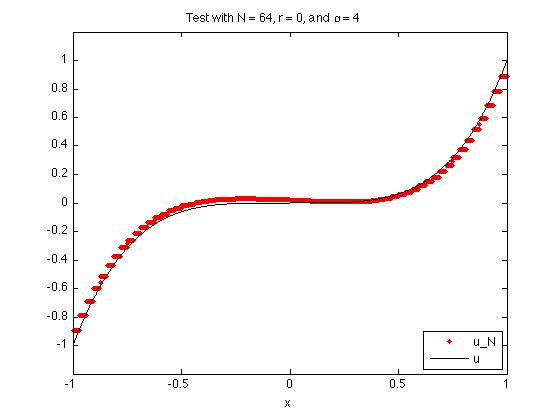}
\includegraphics[scale=0.35]{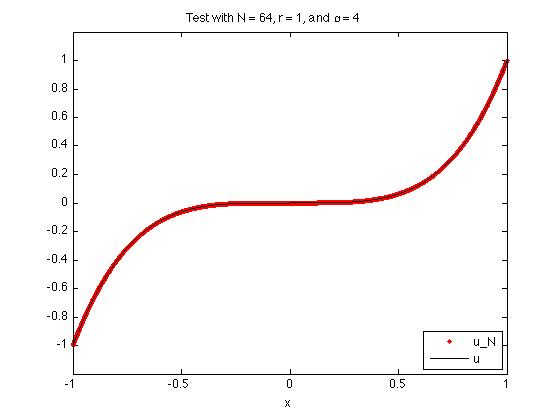}
}
\centerline{
{\tiny
\begin{tabular}{| c | c | c | c c | c c | c c | c c |} 
		\hline
	$r$ & Norm & $h = 1/2$ 
		&\multicolumn{2}{| c |}{$h = 1/4$}
		&\multicolumn{2}{| c |}{$h = 1/8$}
		&\multicolumn{2}{| c |}{$h = 1/16$}
		&\multicolumn{2}{| c |}{$h = 1/32$} \\ 
		\cline{3-11} 
	&   & Error & Error & Order & Error & Order & Error & Order & Error & Order  \\ 
		\hline \cline{1-11}
	0 & $L^{2}$ & 5.0e-01 & 3.4e-01 & 0.53 & 1.3e-01 
	  & 1.37 & 6.1e-02 & 1.12 
	  & 4.4e-02 & 0.49  \\ 
	  & $L^{\infty}$ & 8.5e-01 & 5.7e-01 & 0.57 & 3.5e-01 
	  & 0.72 & 2.0e-01 & 0.79 
	  & 1.1e-01 & 0.86  \\ 
		\hline
	1 & $L^{2}$ & 1.4e-01 & 4.3e-02 & 1.72 & 9.7e-03 
	  & 2.16 & 2.7e-03 & 1.85 
	  & 7.3e-04 & 1.87  \\ 
	  & $L^{\infty}$ & 3.6e-01 & 1.1e-01 & 1.74 & 3.0e-02 
	  & 1.87 & 7.8e-03 & 1.94 
	  & 2.0e-03 & 1.97  \\ 
		\hline
	2 & $L^{2}$ & 2.8e-02 & 3.2e-03 & 3.10 & 4.0e-04 
	  & 3.00 & 5.1e-05 & 2.99 
	  & 6.4e-06 & 2.99  \\ 
	  & $L^{\infty}$ & 4.0e-02 & 5.5e-03 & 2.84 & 7.3e-04 
	  & 2.93 & 9.3e-05 & 2.97 
	  & 1.2e-05 & 2.98  \\ 
		\hline
	3 & $L^{2}$ & 9.4e-03 & 1.3e-03 & 2.91 & 1.6e-04 
	  & 3.01 & 1.9e-05 & 3.01 
	  & 2.4e-06 & 3.01  \\ 
	  & $L^{\infty}$ & 1.1e-02 & 1.5e-03 & 2.91 & 1.9e-04 
	  & 3.02 & 2.3e-05 & 3.01 
	  & 2.9e-06 & 3.01  \\ 
		\hline
\end{tabular}
}
}
\caption{Test 3 with $\alpha = 4$.}
\label{test3_1d}
\end{figure}

\medskip
{\bf Test 3.} 
Consider the stationary Hamilton-Jacobi-Bellman problem with finite dimensional control set
\begin{align*}
	\min_{\theta(x) \in \{1,2 \} } \left\{ 
		- \theta u_{ x x } + u_{x} - u + S(x) \right\}
	& = 0 , \qquad -1 < x < 1 , \\
	u(-1) = -1 , \quad
	u(1) = 1 ,
\end{align*}
where 
\[
	S(x) = \begin{cases}
			-12 x^2 - 4 |x|^3 + x |x|^3 , & x < 0 \\
			24 x^2 - 4 |x|^3 + x |x|^3 , & x \geq 0 .
		\end{cases}
\]
This problem has the exact solution $u(x) = x |x|^3$ corresponding to
$\theta^* (x) = 1$ for $x < 0$ and $\theta^*(x) = 2$ for $x \geq 0$.
Approximating the problem using various order elements, we have the following 
results recorded in Figure~\ref{test3_1d}.  
Again, due to the low regularity of the solution, we expect the rates of convergence
to be affected.  However, we still see increased accuracy for high order elements.

\begin{figure} 
\centerline{
\includegraphics[scale=0.35]{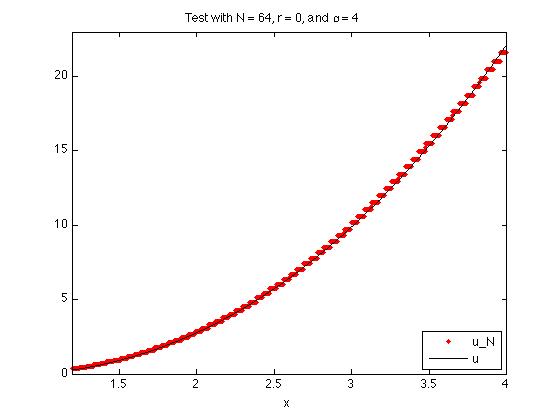}
\includegraphics[scale=0.35]{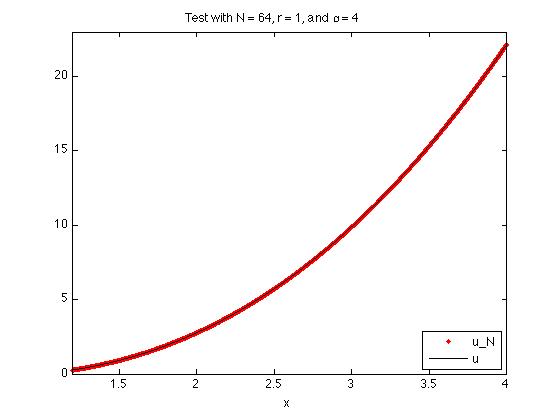}
}
\centerline{
{\tiny
\begin{tabular}{| c | c | c | c c | c c | c c | c c |} 
		\hline
	$r$ & Norm & $h = 2.8/4$ 
		&\multicolumn{2}{| c |}{$h = 2.8/8$}
		&\multicolumn{2}{| c |}{$h = 2.8/16$}
		&\multicolumn{2}{| c |}{$h = 2.8/32$}
		&\multicolumn{2}{| c |}{$h = 2.8/64$} \\ 
		\cline{3-11} 
	&   & Error & Error & Order & Error & Order & Error & Order & Error & Order  \\ 
		\hline \cline{1-11}
	0 & $L^{2}$ & 5.2 & 3.3 & 0.67 & 1.5 
	  & 1.08 & 6.1e-01 & 1.34 
	  & 2.6e-01 & 1.25  \\ 
	  & $L^{\infty}$ & 5.7 & 3.6 & 0.65 & 1.9 
	  & 0.91 & 1.1 & 0.84 
	  & 5.7e-01 & 0.91  \\ 
		\hline
	1 & $L^{2}$ & 2.6e-01 & 8.6e-02 & 1.60 & 2.6e-02 
	  & 1.72 & 7.4e-03 & 1.83 
	  & 2.0e-03 & 1.90  \\ 
	  & $L^{\infty}$ & 3.3e-01 & 1.1e-01 & 1.56 & 3.5e-02 
	  & 1.67 & 1.1e-02 & 1.71 
	  & 3.4e-03 & 1.65  \\ 
		\hline
	2 & $L^{2}$ & 2.6e-03 & 3.9e-04 & 2.77 & 6.6e-05 
	  & 2.55 & 1.4e-05 & 2.26 
	  & 3.2e-06 & 2.09  \\ 
	  & $L^{\infty}$ & 7.3e-03 & 1.0e-03 & 2.85 & 1.4e-04 
	  & 2.91 & 1.9e-05 & 2.81 
	  & 4.1e-06 & 2.25  \\ 
		\hline
	3 & $L^{2}$ & 6.4e-05 & 4.2e-06 & 3.93 & 3.1e-07 
	  & 3.75 & 1.2e-07 & 1.35 
	  & 1.2e-07 & 0.08  \\ 
	  & $L^{\infty}$ & 2.7e-04 & 2.1e-05 & 3.72 & 1.4e-06 
	  & 3.84 & 8.7e-07 & 0.72 
	  & 8.8e-07 & -0.01  \\ 
		\hline
\end{tabular}
}
}
\caption{Test 4 with $\alpha = 4$.}
\label{test4_1d}
\end{figure}

\medskip
{\bf Test 4.} 
Consider the stationary Hamilton-Jacobi-Bellman problem with infinite dimensional control set
\begin{align*}
	\inf_{0 < \theta(x) \leq 1} \left\{ 
		- \theta u_{ x x } + \theta^2 \, x^2 \, u_{x} + \frac{1}{x} u + S(x) \right\}
	& = 0 , & 1.2 < x < 4 , \\
	u(1.2) = 1.44 \ln 1.2 , \quad
	u(4) = 16 \ln 4 ,
\end{align*}
where 
\[
	S(x) = \frac{4 \ln(x)^2 + 12 \ln(x) + 9 - 8 x^4 \ln(x)^2 - 4 x^4 \ln(x)}{4 x^3 \left[ 2 \ln(x) + 1 \right]} .
\]
This problem has the exact solution $u(x) = x^2 \ln x$ corresponding to the control function
$\theta^* (x) = \frac{2 \ln(x) + 3}{2 x^3 \left[ 2 \ln(x) + 1 \right]}$.
Approximating the problem using various order elements,  we obtain the results recorded in Figure~\ref{test4_1d}.
We can see that the approximations appear to reach a maximal level of accuracy of about
5.0e-7 in both $L^2$- and $L^\infty$-norm, which is consistent with the results in Test 
1 corresponding to $r=2$.

\subsection{Parabolic test problems}

We now implement the proposed fully discrete RK4 and trapezoidal LDG methods 
for approximating fully nonlinear parabolic equations of the form \eqref{pde}.
While the above formulation makes no attempt to formally quantify a CFL condition for the RK4 method, 
for the tests we assume a CFL constraint of the form $\Delta t = \kappa_t h^2$, and note that
the constant $\kappa_t$ appears to decrease as the order of the element increases.
Below we implement both the RK4 method and the trapezoidal method for each test problem. 
However, we make no attempt to classify and compare the efficiency of the two methods.
Instead, we focus on testing and demonstrating the usability of both fully discrete schemes
and their promising potentials.  For explicit scheme tests, we record the parameter $\kappa_t$, and for
implicit scheme tests, we record the time step size $\Delta t$.
Note that the row $0^*$ in the figures corresponding to the RK4 method refers to elements with $r=0$ that use 
the standard $L^2$ projection operator in \eqref{rk4}.

\begin{figure} 
\centerline{
\includegraphics[scale=0.35]{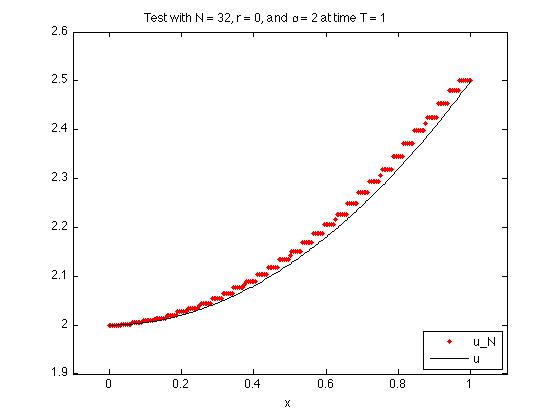}
\hspace{5mm}
\includegraphics[scale=0.35]{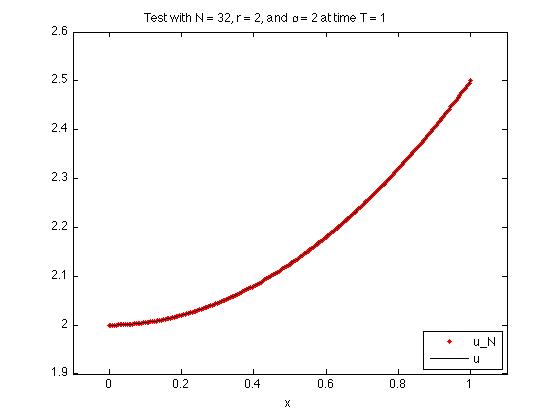}
}
\centerline{
\begin{tabular}{| c | c | c | c c | c c | c c |} 
		\hline
	$r$ & Norm & $h = 1/4$ 
		&\multicolumn{2}{| c |}{$h = 1/8$}
		&\multicolumn{2}{| c |}{$h = 1/16$}
		&\multicolumn{2}{| c |}{$h = 1/32$} \\ 
		\cline{3-9} 
	&   & Error & Error & Order & Error & Order & Error & Order  \\ 
		\hline \cline{1-9}
	0 & $L^{2}$ & 9.9e-02 & 6.4e-02 & 0.64 & 3.6e-02 
	  & 0.81 & 1.9e-02 & 0.92  \\ 
	  & $L^{\infty}$ & 2.2e-01 & 1.4e-01 & 0.64 & 8.0e-02 
	  & 0.81 & 4.3e-02 & 0.89  \\ 
		\hline
	1 & $L^{2}$ & 5.7e-03 & 1.5e-03 & 1.98 & 3.7e-04 
	  & 1.99 & 9.2e-05 & 1.99  \\ 
	  & $L^{\infty}$ & 8.0e-03 & 2.0e-03 & 1.99 & 5.1e-04 
	  & 1.99 & 1.3e-04 & 1.99  \\ 
		\hline
	2 & $L^{2}$ & 2.4e-08 & 2.4e-08 & 0.00 & 2.4e-08 
	  & 0.00 & 2.4e-08 & -0.00  \\ 
	  & $L^{\infty}$ & 3.6e-08 & 3.7e-08 & -0.03 & 3.7e-08 
	  & -0.01 & 3.7e-08 & -0.01  \\ 
		\hline
\end{tabular}
}
\caption{Test 5: Computed solutions at $T = 1$ using
$\kappa_t = 0.001$, $\alpha = 2$.}
\label{exx_2_1_e}
\end{figure}

\begin{figure} 
\centerline{
\includegraphics[scale=0.35]{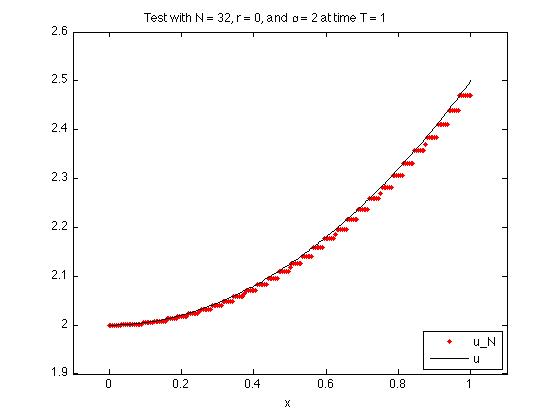}
\hspace{5mm}
\includegraphics[scale=0.35]{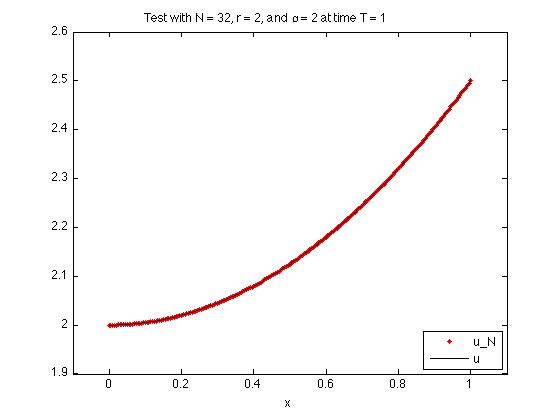}
}
\centerline{
\begin{tabular}{| c | c | c | c c | c c | c c |} 
		\hline
	$r$ & Norm & $h = 1/4$ 
		&\multicolumn{2}{| c |}{$h = 1/8$}
		&\multicolumn{2}{| c |}{$h = 1/16$}
		&\multicolumn{2}{| c |}{$h = 1/32$} \\ 
		\cline{3-9} 
	&   & Error & Error & Order & Error & Order & Error & Order  \\ 
		\hline \cline{1-9}
	0 & $L^{2}$ & 5.9e-02 & 3.4e-02 & 0.78 & 1.9e-02 
	  & 0.84 & 1.0e-02 & 0.91  \\ 
	  & $L^{\infty}$ & 1.9e-01 & 1.1e-01 & 0.79 & 5.9e-02 
	  & 0.90 & 3.0e-02 & 0.95  \\ 
		\hline
	1 & $L^{2}$ & 5.0e-03 & 1.4e-03 & 1.86 & 3.6e-04 
	  & 1.94 & 9.1e-05 & 1.97  \\ 
	  & $L^{\infty}$ & 1.1e-02 & 2.8e-03 & 2.00 & 7.1e-04 
	  & 2.00 & 1.8e-04 & 2.00  \\ 
		\hline
	2 & $L^{2}$ & 1.1e-07 & 1.1e-07 & 0.00 & 1.1e-07 
	  & 0.00 & 1.1e-07 & 0.00  \\ 
	  & $L^{\infty}$ & 1.5e-07 & 1.6e-07 & -0.04 & 1.6e-07 
	  & -0.01 & 1.6e-07 & -0.00  \\ 
		\hline
\end{tabular}
}
\caption{Test 5: Computed solution at $T = 1$ using
$\Delta t = 0.001$ and $\alpha = 2$.}
\label{exx_2_1_i}
\end{figure}

\medskip
{\bf Test 5.} 
Let $\Omega = (0,1)$, $u_a(t) = t^4, u_b = \frac{1}{2} + t^4$, 
and $u_0(x) = \frac{1}{2} x^2$.  
We consider the problem \eqref{pde}, \eqref{ic}, and \eqref{bc_time} with
\[
F(u_{xx}, u_x, u, t, x) = -u_{x x} \, u + \frac{1}{2} x^2 + t^4 - 4 \, t^3 + 1.
\]
It is easy to verify that this problem has a unique classical solution
$u(x,t) = 0.5 \, x^2 + t^4 + 1$. Notice that the PDE has a product 
nonlinearity in the second order derivative. The numerical results
for RK4 are presented in 
Figure~\ref{exx_2_1_e} and the results for the trapezoidal method 
are shown in Figure~\ref{exx_2_1_i}. 
As expected, RK4 recovers the exact solution when $r=2$.

\begin{figure} 
\centerline{
\includegraphics[scale=0.35]{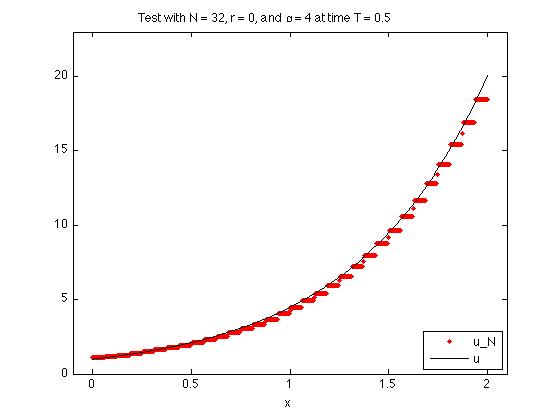}
\hspace{5mm}
\includegraphics[scale=0.35]{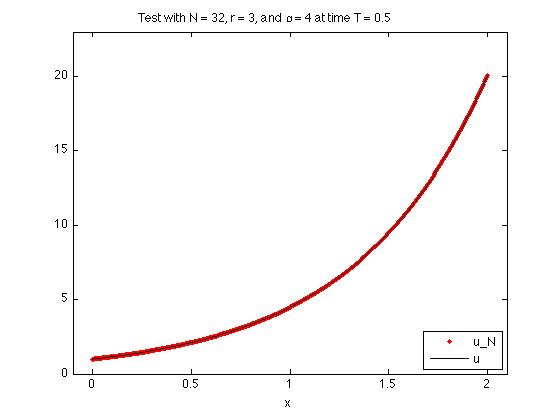}
}
\centerline{
\begin{tabular}{| c | c | c | c c | c c | c c |} 
		\hline
	$r$ & Norm & $h = 2/4$ 
		&\multicolumn{2}{| c |}{$h = 2/8$}
		&\multicolumn{2}{| c |}{$h = 2/16$}
		&\multicolumn{2}{| c |}{$h = 2/32$} \\ 
		\cline{3-9} 
	&   & Error & Error & Order & Error & Order & Error & Order  \\ 
		\hline \cline{1-9}
	0 & $L^{2}$ & 6.9e+00 & 5.7e+00 & 0.28 & 4.1e+00 
	  & 0.47 & 2.6e+00 & 0.65  \\ 
	  & $L^{\infty}$ & 1.0e+01 & 7.9e+00 & 0.40 & 5.6e+00 
	  & 0.50 & 3.7e+00 & 0.60  \\ 
		\hline
	$0^*$ & $L^{2}$ & 2.8e+00 & 1.5e+00 & 0.89 & 8.6e-01 
	  & 0.82 & 5.1e-01 & 0.76  \\ 
	  & $L^{\infty}$ & 7.8e+00 & 5.0e+00 & 0.65 & 2.9e+00 
	  & 0.76 & 1.6e+00 & 0.85  \\ 
		\hline
	1 & $L^{2}$ & 4.8e-01 & 1.2e-01 & 2.05 & 3.0e-02 
	  & 1.93 & 8.2e-03 & 1.88  \\ 
	  & $L^{\infty}$ & 8.4e-01 & 2.3e-01 & 1.87 & 5.8e-02 
	  & 1.99 & 1.5e-02 & 1.92  \\ 
		\hline
	2 & $L^{2}$ & 3.7e-02 & 7.8e-03 & 2.23 & 1.9e-03 
	  & 2.03 & 4.8e-04 & 2.00  \\ 
	  & $L^{\infty}$ & 5.9e-02 & 1.0e-02 & 2.53 & 2.2e-03 
	  & 2.19 & 5.2e-04 & 2.11  \\ 
		\hline
	3 & $L^{2}$ & 1.1e-03 & 7.4e-05 & 3.95 & 4.7e-06 
	  & 3.99 & 3.01e-07 & 3.96  \\ 
	  & $L^{\infty}$ & 2.5e-03 & 1.6e-04 & 3.93 & 1.1e-05 
	  & 3.86 & 7.66e-07 & 3.84  \\ 
		\hline
\end{tabular}
}
\caption{Test 6: Computed solutions at time $T = 0.5$ using
$\kappa_t = 0.005, 0.001, 0.0005, 0.0001$ for $r=0, 1, 2, 3$, respectively, and $\alpha = 4$.
Left figure uses the standard $L^2$ projection operator.}
\label{exx_2_2_e}
\end{figure}

\begin{figure}
\centerline{
\includegraphics[scale=0.35]{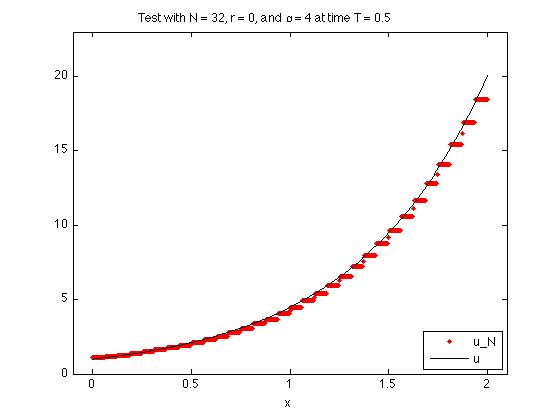}
\hspace{5mm}
\includegraphics[scale=0.35]{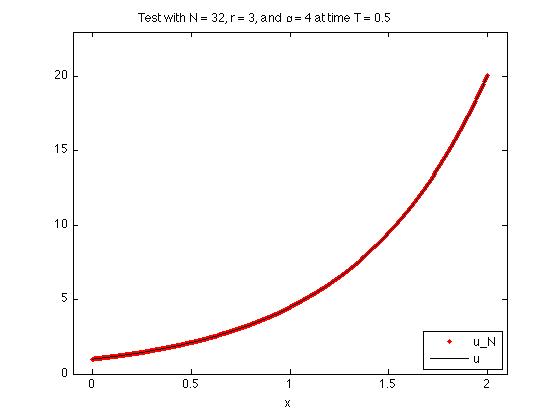}
}
\centerline{
\begin{tabular}{| c | c | c | c c | c c | c c |} 
		\hline
	$r$ & Norm & $h = 2/4$ 
		&\multicolumn{2}{| c |}{$h = 2/8$}
		&\multicolumn{2}{| c |}{$h = 2/16$}
		&\multicolumn{2}{| c |}{$h = 2/32$} \\ 
		\cline{3-9} 
	&   & Error & Error & Order & Error & Order & Error & Order  \\ 
		\hline \cline{1-9}
	0 & $L^{2}$ & 2.8e+00 & 1.5e+00 & 0.89 & 8.6e-01 
	  & 0.82 & 5.1e-01 & 0.77  \\ 
	  & $L^{\infty}$ & 7.8e+00 & 5.0e+00 & 0.65 & 2.9e+00 
	  & 0.76 & 1.6e+00 & 0.85  \\ 
		\hline
	1 & $L^{2}$ & 3.8e-01 & 1.3e-01 & 1.62 & 4.5e-02 
	  & 1.49 & 1.5e-02 & 1.60  \\ 
	  & $L^{\infty}$ & 1.3e+00 & 4.0e-01 & 1.74 & 1.1e-01 
	  & 1.85 & 3.0e-02 & 1.91  \\ 
		\hline
	2 & $L^{2}$ & 2.7e-02 & 6.7e-03 & 2.04 & 1.9e-03 
	  & 1.82 & 5.3e-04 & 1.85  \\ 
	  & $L^{\infty}$ & 1.0e-01 & 1.5e-02 & 2.77 & 2.1e-03 
	  & 2.80 & 5.4e-04 & 1.99  \\ 
		\hline
	3 & $L^{2}$ & 1.1e-03 & 7.2e-05 & 3.89 & 1.3e-05 
	  & 2.46 & 1.2e-05 & 0.09  \\ 
	  & $L^{\infty}$ & 5.5e-03 & 4.0e-04 & 3.76 & 2.7e-05 
	  & 3.92 & 1.3e-05 & 1.05  \\ 
		\hline
\end{tabular}
}
\caption{Test 6: Computed solutions at time $T = 0.5$ using
$\Delta t = 0.005$ and $\alpha = 4$.}
\label{exx_2_2_i}
\end{figure}

\medskip
{\bf Test 6.} 
Let $\Omega = (0,2)$, $u_a(t) = 1$, $u_b = e^{2 (t+1)}$, and $u_0(x) = e^x$.  
We consider the problem \eqref{pde}, \eqref{ic}, and \eqref{bc_time} with
\[
F(u_{x x}, u_x, u, t, x) = - u_x \ln \bigl( u_{x x} + 1 \bigr) + S(x,t) , 
\]
and
\[
S(x,t) = e^{(t+1) x} \Big(x - (t+1) \ln \bigl((t+1)^2 e^{(t+1) x}+1 \bigr) \Big).
\]
It is easy to verify that this problem has a unique classical
solution $u(x,t) = e^{(t+1) x}$. Notice that this problem is nonlinear in 
both $u_{x x}$ and $u_x$. Furthermore, the exact solution $u$ cannot be 
factored into the form $u(x,t) = G(t) \, Y (x)$ for some functions $G$ and $Y$.
Results for RK4 are given in Figure~\ref{exx_2_2_e}, and results for the 
trapezoidal method are shown in Figure~\ref{exx_2_2_i}.
We note that RK4 was unstable without using the very restrictive values for $\kappa_t$
recorded in Figure~\ref{exx_2_2_e}.  
However, for RK4, we observe optimal rates of convergence in the spacial variable  
while the rates for the trapezoidal method appear to be limited by the lower rate of 
convergence for the time-stepping scheme.

\begin{figure} 
\centerline{
\includegraphics[scale=0.35]{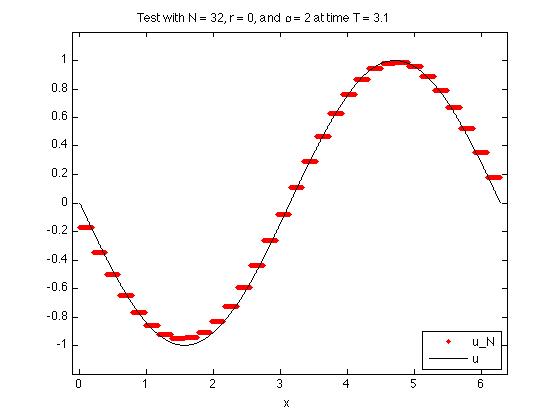}
\hspace{5mm}
\includegraphics[scale=0.35]{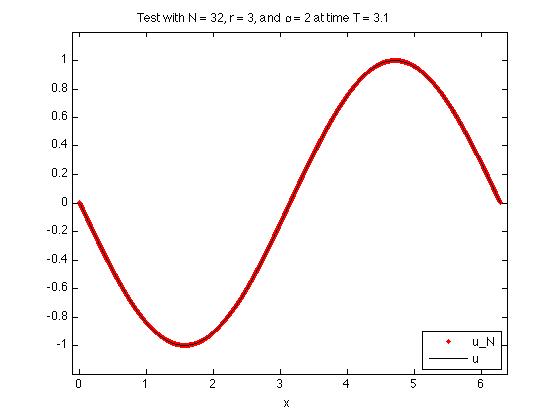}
}
\centerline{
\begin{tabular}{| c | c | c | c c | c c | c c |} 
		\hline
	$r$ & Norm & $h = \pi/2$ 
		&\multicolumn{2}{| c |}{$h = \pi/4$}
		&\multicolumn{2}{| c |}{$h = \pi/8$}
		&\multicolumn{2}{| c |}{$h = \pi/16$} \\ 
		\cline{3-9} 
	&   & Error & Error & Order & Error & Order & Error & Order  \\ 
		\hline \cline{1-9}
	0 & $L^{2}$ & 1.5e+00 & 1.3e+00 & 0.24 & 7.1e-01 
	  & 0.82 & 3.3e-01 & 1.12  \\ 
	  & $L^{\infty}$ & 9.6e-01 & 7.7e-01 & 0.31 & 4.9e-01 
	  & 0.65 & 2.9e-01 & 0.78  \\ 
		\hline
	$0^*$ & $L^{2}$ & 1.2e+00 & 6.9e-01 & 0.78 & 3.0e-01 
	  & 1.19 & 1.4e-01 & 1.16  \\ 
	  & $L^{\infty}$ & 9.2e-01 & 5.8e-01 & 0.67 & 3.2e-01 
	  & 0.85 & 1.8e-01 & 0.83  \\ 
		\hline
	1 & $L^{2}$ & 2.7e-01 & 7.6e-02 & 1.81 & 2.0e-02 
	  & 1.94 & 5.0e-03 & 1.99  \\ 
	  & $L^{\infty}$ & 1.9e-01 & 6.6e-02 & 1.52 & 1.7e-02 
	  & 1.97 & 4.2e-03 & 2.00  \\ 
		\hline
	2 & $L^{2}$ & 7.2e-02 & 1.8e-02 & 1.98 & 4.5e-03 
	  & 2.02 & 1.1e-03 & 2.01  \\ 
	  & $L^{\infty}$ & 6.8e-02 & 1.8e-02 & 1.93 & 4.3e-03 
	  & 2.04 & 1.1e-03 & 2.03  \\ 
		\hline
	3 & $L^{2}$ & 8.3e-03 & 5.7e-04 & 3.87 & 3.6e-05 
	  & 3.97 & 2.2e-06 & 4.02  \\ 
	  & $L^{\infty}$ & 7.6e-03 & 5.1e-04 & 3.92 & 3.2e-05 
	  & 4.00 & 1.9e-06 & 4.04  \\ 
		\hline
\end{tabular}
}
\caption{Test 7: Computed solutions at time $T = 3.10$ using
$\kappa_t = 0.05, 0.005, 0.001, 0.0005$ for $r=0, 1, 2, 3$, respectively, and $\alpha = 2$.  
Left plot uses the standard $L^2$ projection operator.
}
\label{exx_2_3_e}
\end{figure}

\begin{figure}
\centerline{
\includegraphics[scale=0.35]{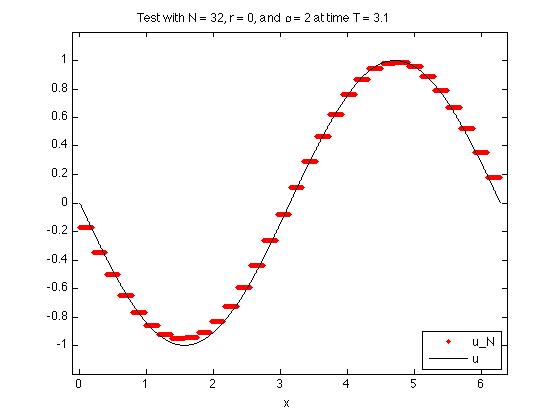}
\hspace{5mm}
\includegraphics[scale=0.35]{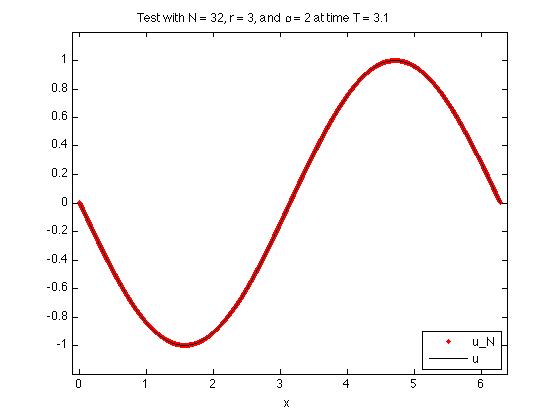}
}
\centerline{
\begin{tabular}{| c | c | c | c c | c c | c c |} 
		\hline
	$r$ & Norm & $h = \pi/2$ 
		&\multicolumn{2}{| c |}{$h = \pi/4$}
		&\multicolumn{2}{| c |}{$h = \pi/8$}
		&\multicolumn{2}{| c |}{$h = \pi/16$} \\ 
		\cline{3-9} 
	&   & Error & Error & Order & Error & Order & Error & Order  \\ 
		\hline \cline{1-9}
	0 & $L^{2}$ & 1.2e+00 & 6.9e-01 & 0.78 & 3.0e-01 
	  & 1.19 & 1.4e-01 & 1.16  \\ 
	  & $L^{\infty}$ & 9.2e-01 & 5.8e-01 & 0.67 & 3.2e-01 
	  & 0.85 & 1.8e-01 & 0.83  \\ 
		\hline
	1 & $L^{2}$ & 2.3e-01 & 7.5e-02 & 1.63 & 2.0e-02 
	  & 1.89 & 7.9e-03 & 1.36  \\ 
	  & $L^{\infty}$ & 2.1e-01 & 6.6e-02 & 1.66 & 1.7e-02 
	  & 1.95 & 5.8e-03 & 1.56  \\ 
		\hline
	2 & $L^{2}$ & 6.9e-02 & 1.8e-02 & 1.93 & 4.7e-03 
	  & 1.95 & 2.5e-03 & 0.90  \\ 
	  & $L^{\infty}$ & 6.5e-02 & 1.8e-02 & 1.87 & 4.5e-03 
	  & 2.00 & 1.9e-03 & 1.21  \\ 
		\hline
	3 & $L^{2}$ & 8.1e-03 & 5.9e-04 & 3.79 & 1.1e-04 
	  & 2.42 & 1.1e-04 & 0.03  \\ 
	  & $L^{\infty}$ & 7.5e-03 & 5.2e-04 & 3.87 & 7.6e-05 
	  & 2.77 & 7.3e-05 & 0.06  \\ 
		\hline
\end{tabular}
}
\caption{Test 7: Computed solutions at time $T = 3.10$ using
$\Delta t = 0.031$ and $\alpha = 2$.}
\label{exx_2_3_i}
\end{figure}

\medskip
{\bf Test 7.} 
Let $\Omega = (0, 2 \pi)$, $u_a(t) = 0$, $u_b = 0$, and $u_0(x) = \sin (x)$.  
We consider the problem \eqref{pde}, \eqref{ic}, and \eqref{bc_time} with
\[
F(u_{x x}, u_x, u, t, x) = - \min_{\theta (t, x) \in \{1,2 \}} 
\Big\{A_\theta\, u_{x x}-c(x,t) \cos (t) \, \sin (x) - \sin(t) \, \sin(x) \Big\} , 
\]
where $A_1 = 1, A_2 = \frac12$, and
\[
c(x,t) = \begin{cases}
1, & \text{if } 
0 < t \leq \frac{\pi}2 \mbox{ and } 0 < x \leq \pi \text{ or } 
\frac{\pi}2 < t \leq \pi \text{ and } \pi < x < 2 \pi, \\
\frac12, & \mbox{otherwise}. 
\end{cases}
\]
It is easy to check that this problem has a unique classical solution
$u (x,t) = \cos(t) \sin(x)$.  Notice that this problem has a finite 
dimensional control parameter set, and the optimal control is given by
\[
\theta^* (t, x) = \begin{cases}
1, & \mbox{if } c(x,t) = 1 , \\
2 , & \mbox{if } c(x,t) = \frac12 .
\end{cases} 
\]
The numerical results are reported in Figure~\ref{exx_2_3_e} for RK4 
and in Figure~\ref{exx_2_3_i} for the trapezoidal method.

\begin{figure} 
\centerline{
\includegraphics[scale=0.35]{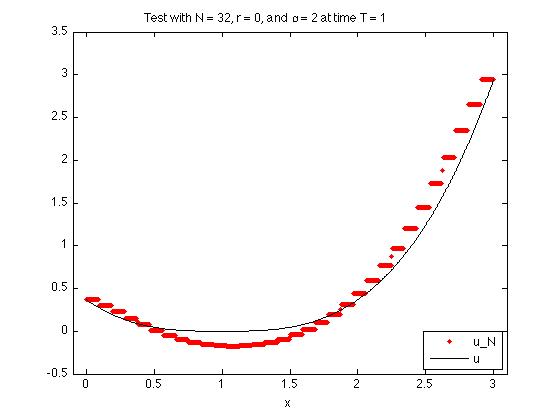}
\hspace{5mm}
\includegraphics[scale=0.35]{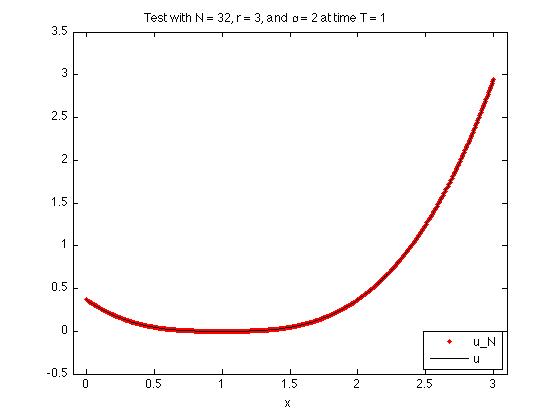}
}
\centerline{
\begin{tabular}{| c | c | c | c c | c c | c c |} 
		\hline
	$r$ & Norm & $h = 3/4$ 
		&\multicolumn{2}{| c |}{$h = 3/8$}
		&\multicolumn{2}{| c |}{$h = 3/16$}
		&\multicolumn{2}{| c |}{$h = 3/32$} \\ 
		\cline{3-9} 
	&   & Error & Error & Order & Error & Order & Error & Order  \\ 
		\hline \cline{1-9}
	0 & $L^{2}$ & 2.1e+00 & 1.4e+00 & 0.51 & 5.3e-01 
	  & 1.44 & 2.9e-01 & 0.88  \\ 
	  & $L^{\infty}$ & 2.1e+00 & 1.5e+00 & 0.44 & 8.5e-01 
	  & 0.86 & 4.9e-01 & 0.81  \\ 
		\hline
	$0^*$ & $L^{2}$ & 8.6e-01 & 8.8e-01 & -0.02 & 5.9e-01 
	  & 0.58 & 2.9e-01 & 1.00  \\ 
	  & $L^{\infty}$ & 1.8e+00 & 1.3e+00 & 0.53 & 7.3e-01 
	  & 0.80 & 3.9e-01 & 0.92  \\ 
		\hline
	1 & $L^{2}$ & 2.0e-01 & 7.3e-02 & 1.45 & 1.3e-02 
	  & 2.44 & 3.2e-03 & 2.06  \\ 
	  & $L^{\infty}$ & 3.9e-01 & 1.0e-01 & 1.92 & 2.7e-02 
	  & 1.95 & 6.8e-03 & 1.98  \\ 
		\hline
	2 & $L^{2}$ & 1.1e-01 & 2.0e-02 & 2.49 & 5.2e-03 
	  & 1.98 & 1.3e-03 & 2.02  \\ 
	  & $L^{\infty}$ & 1.1e-01 & 2.0e-02 & 2.42 & 5.2e-03 
	  & 1.99 & 1.3e-03 & 2.00  \\ 
		\hline
	3 & $L^{2}$ & 3.0e-02 & 6.8e-03 & 2.13 & 1.7e-03 
	  & 2.00 & 4.2e-04 & 2.03  \\ 
	  & $L^{\infty}$ & 2.9e-02 & 6.9e-03 & 2.08 & 1.7e-03 
	  & 2.02 & 4.2e-04 & 2.02  \\ 
		\hline
\end{tabular}
}
\caption{Test 8: Computed solutions at time $T = 1$ using
$\kappa_t = 0.05, 0.005, 0.001, 0.0005$ for $r=0, 1, 2, 3$, respectively, and $\alpha = 2$.
}
\label{exx_2_4_e}
\end{figure}

\begin{figure}
\centerline{
\includegraphics[scale=0.35]{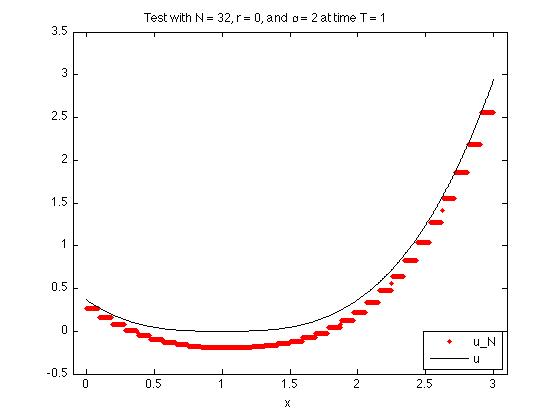}
\hspace{5mm}
\includegraphics[scale=0.35]{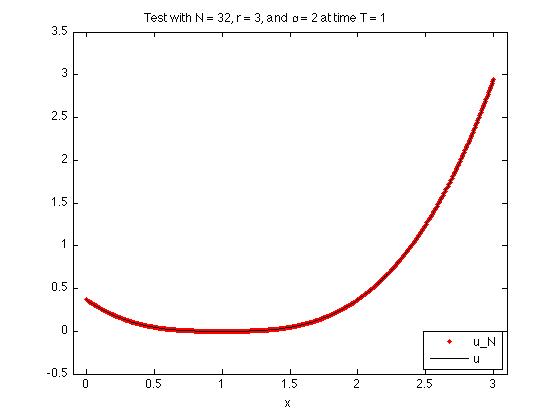}
}
\centerline{
\begin{tabular}{| c | c | c | c c | c c | c c |} 
		\hline
	$r$ & Norm & $h = 3/4$ 
		&\multicolumn{2}{| c |}{$h = 3/8$}
		&\multicolumn{2}{| c |}{$h = 3/16$}
		&\multicolumn{2}{| c |}{$h = 3/32$} \\ 
		\cline{3-9} 
	&   & Error & Error & Order & Error & Order & Error & Order  \\ 
		\hline \cline{1-9}
	0 & $L^{2}$ & 8.6e-01 & 8.8e-01 & -0.02 & 5.9e-01 
	  & 0.58 & 2.9e-01 & 1.00  \\ 
	  & $L^{\infty}$ & 1.8e+00 & 1.3e+00 & 0.53 & 7.3e-01 
	  & 0.80 & 3.9e-01 & 0.92  \\ 
		\hline
	1 & $L^{2}$ & 2.0e-01 & 7.3e-02 & 1.45 & 1.3e-02 
	  & 2.44 & 3.2e-03 & 2.05  \\ 
	  & $L^{\infty}$ & 3.9e-01 & 1.0e-01 & 1.92 & 2.7e-02 
	  & 1.95 & 6.8e-03 & 1.98  \\ 
		\hline
	2 & $L^{2}$ & 1.1e-01 & 2.0e-02 & 2.49 & 5.2e-03 
	  & 1.98 & 1.3e-03 & 2.02  \\ 
	  & $L^{\infty}$ & 1.1e-01 & 2.0e-02 & 2.42 & 5.2e-03 
	  & 1.99 & 1.3e-03 & 2.00  \\ 
		\hline
	3 & $L^{2}$ & 3.0e-02 & 6.8e-03 & 2.13 & 1.7e-03 
	  & 2.00 & 4.2e-04 & 2.03  \\ 
	  & $L^{\infty}$ & 2.9e-02 & 6.9e-03 & 2.08 & 1.7e-03 
	  & 2.02 & 4.2e-04 & 2.02  \\ 
		\hline
\end{tabular}
}
\caption{Test 8: Computed solutions at time $T = 1$ using
$\Delta t = 0.001$ and $\alpha = 2$.}
\label{exx_2_4_i}
\end{figure}

\medskip
{\bf Test 8.} 
Let $\Omega = (0,3)$, $u_a(t) = e^{-t}$, $u_b = 8 \, e^{-t}$, and $u_0(x) = | x - 1 |^3$.  
We consider the problem \eqref{pde}, \eqref{ic}, and \eqref{bc_time} with
\[
	F(u_{x x}, u_x, u, t, x) = - \inf_{-1 \leq \theta(t,x) \leq 1} \Big\{ | x-1 | u_{x x} + \theta u_x \Big\} 
		- | x-1 |^2 \, \left( | x-1 | - 3 \right) \, e^{-t} , 
\]
It is easy to verify that the problem has the exact solution
$u(t,x) =  | x - 1 |^3 \, e^{-t}$.
Notice that the above operator is not elliptic for $x = 1$.  
Also, this problem has a bang-bang type control with the optimal control given by
\[
	\theta^* (t,x) = \begin{cases}
		1, & \text{if } x < 1, \\
		-1 , & \text{if } x > 1. 
	\end{cases}
\]
We can see from the results for the RK4 method in Figure~\ref{exx_2_4_e} and the results for
the trapezoidal method in Figure~\ref{exx_2_4_i} that the spatial rates of convergence
appear to be limited by the low regularity of the solution, while the accuracy appears to 
increase with respect to the element degree.


\subsection{The role of the numerical moment} \label{alpha_tests}

In this section, we focus on understanding the role of the numerical moment.
We first give an interpretation of the numerical moment,  and then we 
explore the utility of the numerical moment.
The role of the numerical moment can heuristically be understood as follows
when the numerical moment is rewritten in the form
\[
h^2 \, \alpha \, \Bigl(\frac{ p_{1h} - p_{2h} - p_{3h} + p_{4h}}{h^2}\Bigr).
\]
Letting $r = 0$, we see that $\frac{p_{1h} - p_{2h} - p_{3h} + p_{4h}}{h^2}$ is
an $\mathcal{O}( h^2)$ approximation to $u_{x x x x}$.
Then, we can see that the numerical moment acts as a centered difference
approximation for $\Delta^2 u$ multiplied by a factor that tends to zero with rate
$\mathcal{O}( h^2)$. Thus, at the PDE level, we are in essence approximating 
the fully nonlinear second order elliptic operator
\[
F \left( u_{x x}, u_x, u, x \right)
\]
by the quasilinear fourth order operator $\widetilde{F}_\rho$, where
\[
\widetilde{F}_\rho \left( u_{x x x x}, u_{x x}, u_x, u, x \right) 
=\rho \, u_{x x x x} + F \left( u_{x x}, u_x, u, x \right) .
\]  
In the limit as $\rho \to 0$, we heuristically expect that the solution 
of the fourth order problem converges to the unique viscosity solution 
of the second order problem.  Using a convergent family of fourth order 
quasilinear PDEs to approximate a fully nonlinear second order PDE 
has previously been considered for PDEs such as the Monge-Amp\`{e}re equation, 
the prescribed Gauss curvature equation, the infinity-Laplace equation, 
and linear second order equations of non-divergence form.  The method is known
as the vanishing moment method. We refer the reader to 
\cite{Feng_Neilan11,FGN12} for a detailed exposition.
  
We now show that the proposed schemes using a numerical moment can
eliminate the numerical artifacts
that arise as consequences from using a standard discretization, and in certain cases when the
numerical artifacts are not fully eliminated, the algebraic system has enough structure
to design solvers that emphasize the monotonicity in $\frac{p_2+p_3}{2}$ 
and are consistent in searching for solutions at which the nonlinear PDE problem is elliptic. 
We again consider the Monge-Amp\`{e}re type problem from Test 1 in section~\ref{ell_tests_1d}.
The given problem has two classical PDE solutions $u^+$ and $u^-$.  
However, there exists infinitely many $C^1$ functions that satisfy the PDE and boundary
conditions almost everywhere, as seen by $\mu$ in \eqref{art}.
These almost everywhere solutions will correspond to what we call numerical artifacts 
in that algebraic solutions for a given discretization may correspond to these functions.
It is well known that using a standard discretization scheme for the Monge-Amp\`{e}re problem
can yield multiple solutions, many of which are numerical artifacts that do not correspond to 
PDE solutions \cite{FGN12}. For example, let $\mu \in  H^2(0,1) \setminus C^2(0,1)$ be defined by
\begin{equation} \label{art}
	\mu(x) = \begin{cases}
		\frac{1}{2} x^2 + \frac{1}{4} x,  & \text{for } x < 0.5 , \\
		- \frac{1}{2} x^2 + \frac{5}{4} x - \frac{1}{4},  & \text{for } x \geq 0.5 .
		\end{cases}
\end{equation}
Furthermore, suppose $x_j = 0.5$ for some $j = 2, 3, \ldots, J-1$.
Then, when using a standard central difference discretization, $\mu$ corresponds to a numerical 
solution, yielding a numerical artifact.

We now consider our discretization that uses a numerical moment.  When $\alpha=0$, we have
no numerical moment.  As a result, we have numerical artifacts as in the standard central difference 
discretization case. Suppose $r=0$.  Then, for $\alpha \neq 0$, inspection of \eqref{p1h_r0}--\eqref{p4h_r0}
yields the fact that $p_{2h}$ cannot jump from a value of $1$ to a value of $-1$ when 
going across $x_j=0.5$.  Thus, the numerical moment penalizes discontinuities in $p_{jh}$, 
$j=1,2,3,4$,  and, as a result, the numerical moment eliminates numerical artifacts such as $\mu$.
Similarly, for $r=1$, we can see that $\mu$ does not correspond to a numerical solution.
However, in this case, the algebraic system does have a small residual that may trap solvers
such as {\it fsolve}.  
Thus, for $r=0$ and $r=1$, the numerical moment penalizes differences in $p_{jh}$, $j=1,2,3,4$, 
that arise from discontinuities in $u_h$, $q_{1h}$, and $q_{2h}$. Hence,  it eliminates
numerical artifacts such as $\mu$.

We now consider $r \geq 2$, in which case $\mu \in V^h$.  Furthermore, since $\mu \in C^1$, 
we will end up with $u_h = \mu$, $q_{1h} = q_{2h}$, and $p_{1h} = p_{jh}$ for $j=2,3,4$
being a numeric solution, where
\[
q_{1h}(x) = \begin{cases}
	x+\frac14, & \text{for } x < 0.5, \\
	-x + \frac54, & \text{for } x > 0.5 , 
	\end{cases}
\qquad \text{and} \qquad
p_{1h}(x) = \begin{cases}
	1 , & \text{for } x < 0.5, \\
	-1, & \text{for } x > 0.5 . 
	\end{cases}
\]
Thus, by the equalities of $p_{jh}$, $j=1,2,3,4$, the numerical moment is always zero and
we do have numerical artifacts. 
These equalities are a consequence of the continuity of $u_h$, $q_{1h}$, and $q_{2h}$.
With the extra degrees of freedom for $r\geq2$, we allow $C^1$ to be embedded into our
approximation space $V^h$, thus creating possible solutions with a zero valued numerical
moment.  The numerical moment acts as a penalty term for differences in $p_{jh}$, $j=1,2,3,4$, 
which are a consequence of discontinuities in $q_{1h}$ and $q_{2h}$ that naturally arise
for nontrivial functions when $r=0$ or $r=1$.

Even with the possible presence of numerical artifacts for the above discretization when 
$r \geq 2$, the numerical moment can be exploited at the solver level.  We now present a 
splitting algorithm for solving the resulting nonlinear algebraic system that uses the 
numerical moment to strongly emphasize the fact that the viscosity solution of the PDE should
preserve the monotonicity required by the definition of ellipticity. 

\medskip
\begin{solver} \label{solver_alg}
\noindent
\begin{enumerate}
\item[{\rm (1)}] Pick an initial guess for $u_h$. 
\item[{\rm (2)}] Form initial guesses for $q_{1h}$, $q_{2h}$, $p_{1h}$, $p_{2h}$, $p_{3h}$, 
	and $p_{4h}$ using equations \eqref{q_bil} and \eqref{p_bil}.
\item[{\rm (3)}] Solve equation \eqref{pde_ell_weak} for $\frac{p_{2h}+p_{3h}}{2}$. 
\item[{\rm (4)}] Solve equation \eqref{q_bil} for $i=1,2$ and the equation formed
	by averaging \eqref{p_bil} for $j=2,3$ for $u_h$, $q_{1h}$, and $q_{2h}$. 
\item[{\rm (5)}] Solve equation \eqref{p_bil} for $j=1, 2, 3, 4$ for 
	$p_{1h}$, $p_{2h}$, $p_{3h}$, and $p_{4h}$. 
\item[{\rm (6)}] Repeat Steps 3 - 5 until the change in $\frac{p_{2h}+p_{3h}}{2}$ is sufficiently small.
\end{enumerate}
\end{solver}

\begin{figure}
\centerline{
\includegraphics[scale=0.3]{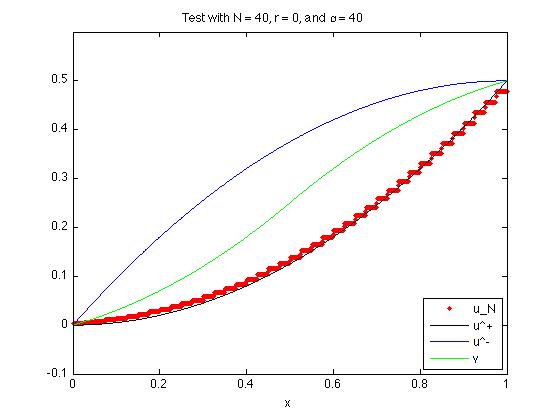}
\hspace{5mm}
\includegraphics[scale=0.3]{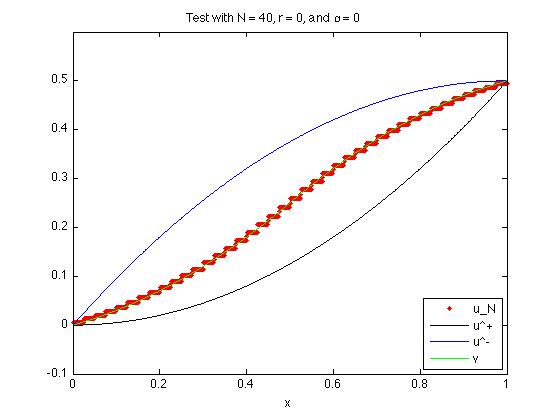}
}
\caption{Left: $\alpha = 40$, $h = 1/40$, and $r=0$. 
Right: $\alpha = 0$, $h = 1/40$, and $r=0$}
\label{alpha_r0}
\end{figure}

\begin{figure}
\centerline{
\includegraphics[scale=0.3]{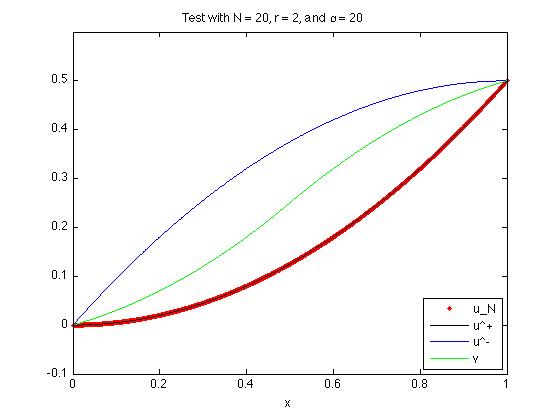}
\hspace{5mm}
\includegraphics[scale=0.3]{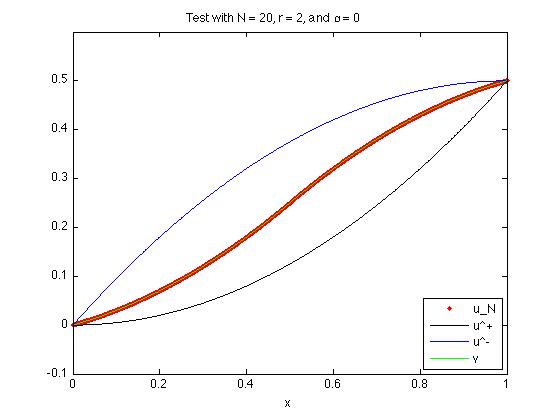}
}
\caption{Left: $\alpha = 20$, $h = 1/20$, and $r=2$. 
Right: $\alpha = 0$, $h = 1/20$, and $r=2$}
\label{alpha_r2}
\end{figure}

\begin{figure}
\centerline{
\includegraphics[scale=0.3]{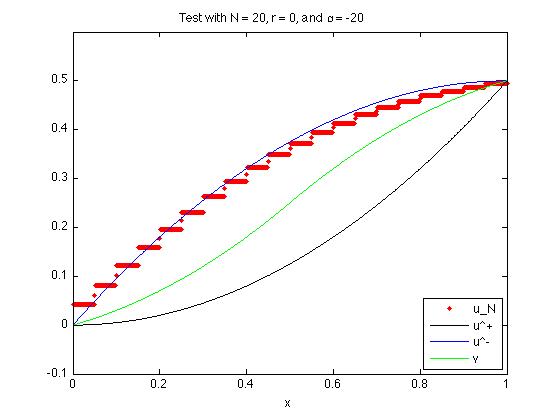}
\hspace{5mm}
\includegraphics[scale=0.3]{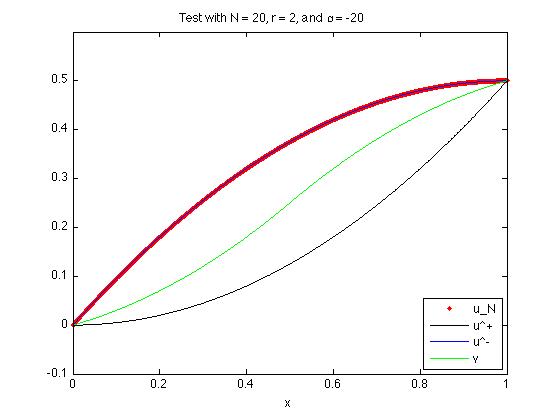}
}
\caption{Left: $\alpha = -40$, $h = 1/40$, and $r=0$. 
Right: $\alpha = -20$, $h = 1/20$, and $r=2$}
\label{alpha_neg}
\end{figure}

\smallskip
For the next numerical tests, we will show that using Algorithm~\ref{solver_alg} with a 
sufficiently large coefficient for the numerical moment destabilizes numerical artifacts
such as $\mu$ and steers the approximation towards the viscosity solution of the PDE.
Let $\ou(x) = \frac{x}2$.  Then, $\ou$ is the secant line formed by the boundary data for the
given boundary value problem.
We now approximate the solution of the Monge-Amp\`{e}re type problem from 
Test 1 in section~\ref{ell_tests_1d}
by using 100 iterations of Algorithm~\ref{solver_alg} followed by using {\it fsolve} on the
full system to solve the global discretization given by 
\eqref{pde_ell_weak}, \eqref{q_bil}, and \eqref{p_bil}.
We take the initial guess to be
\[
	u_{h}^{(0)} = \frac{3}{4} \mu + \frac{1}{4} \ou , 
\]
where, for $r=0$, $u_h^{0}$ is first projected into $V^h$.
From Figure~\ref{alpha_r0}, we see that the numerical moment drives the solution towards 
the viscosity solution $u^+$ when $r = 0$ and $\alpha$ is positive. From Figure~\ref{alpha_r2}, 
we see that the numerical moment also drives the solution towards the viscosity
solution $u^+$ when $r = 2$ and $\alpha$ is positive, despite the presence of numerical artifacts.  
From Figure~\ref{alpha_neg}, we see that the moment drives the solution towards the viscosity
solution of $F(u_{x x}, u_x , u , x) := u_{x x}^2 - 1$ given by $u^-$ for $r = 0$ and $r=2$ when
$\alpha$ is chosen to be negative.
In each figure, the middle graph corresponds to $\mu$.
Clearly, we recover the numerical artifact corresponding to $\mu$ when $\alpha=0$.
Thus, the numerical moment plays an essential role in either eliminating numerical artifacts
at the discretization level or handling numerical artifacts at the solver level.

We make one final note about using the iterative solver given by Algorithm~\ref{solver_alg}.
Note that using {\it fsolve} to solve the full system with the initial guess given by $u_h^{(0)}$ 
resulted in either not finding a root for many tests ($r=1$) or converging to a numerical artifact 
with a discontinuous second order derivative at another node in the mesh ($r=2$).
In order to use {\it fsolve} for the given test problem, the initial guess should either be 
restricted to the class of functions where $p_{2h}$ and $p_{3h}$ preserve the ellipticity of the 
nonlinear operator, the initial guess should be preconditioned by first using {\it fsolve} with
$r = 0,1$, or the initial guess should be preconditioned using Algorithm~\ref{solver_alg}.
When using $r=0$ and a non-ellipticity-preserving initial guess, 
solving the full system of equations with {\it fsolve} still has the potential to converge to 
$u^-$ even for $\alpha > 0$.
The strength of Algorithm~\ref{solver_alg} is that it strongly enforces the requirement that $\hF$ is 
monotone decreasing in $p_2$ and $p_3$ over each iteration.  Thus, a sufficiently large
value for $\alpha$ drives the approximation towards the class of ellipticity-preserving functions if the
algorithm converges.

\bibliographystyle{plain}

\end{document}